\documentclass[11pt]{article}

\usepackage{amsmath,amssymb,latexsym, amsfonts, amscd, amsthm}
\usepackage[small,nohug,heads=littlevee]{diagrams}
\diagramstyle[labelstyle=\scriptstyle]

\newcommand{\reg}{^{\mathrm{reg}}}
\newcommand{\ab}{_{\mathrm{ab}}}

\DeclareMathOperator{\Hom}{Hom}

\theoremstyle{plain}

\newtheorem{theorem}{Theorem}[section]
\newtheorem{conjecture}[theorem]{Conjecture}
\newtheorem{proposition}[theorem]{Proposition}

\newtheorem{lemma}[theorem]{Lemma}

\def\Z{\mathbb{Z}}

\def\C{\mathbb{C}}

\def\R{\mathbb{R}}
\def\Ad{\mathrm{Ad}}

\DeclareMathOperator{\Gal}{Gal}
\DeclareMathOperator{\Spin}{Spin}
\DeclareMathOperator{\im}{im}

\theoremstyle{definition}
\newtheorem{definition}[theorem]{Definition}
\newtheorem{remark}[theorem]{Remark}

\usepackage[pdftex,plainpages=false]{hyperref}

\begin{document}
\title{A Real Groups Construction of the Tame Local Langlands Correspondence for $PGSp(4,F)$}

\author{Moshe Adrian \\ Joshua Lansky}

\maketitle

\begin{abstract}
In this paper, we continue the work in \cite{adrian1} and give a new construction of the tame local Langlands correspondence for $PGSp(4,F)$, where $F$ is a $p$-adic field, that is analogous to the construction of the local Langlands correspondence for real groups.
\end{abstract}

\section{Introduction}
In this paper, we give a new construction of the tame local Langlands correspondence for $PGSp(4,F)$, where $F$ is a non-Archimedean local field of characteristic zero, using character theory and ideas from the theory of real reductive groups.  We assume that the residual characteristic of $F$ is greater than 46 (see Remark \ref{nonvanishing1}).  We continue the program initiated in \cite{adrian1}, where Adrian gave a new realization of the local Langlands correspondence for $GL(\ell,F)$, $\ell$ a prime.

There has been a significant amount of progress in the local Langlands correspondence in recent years.  In the recent work of \cite{bushnellhenniart}, \cite{debackerreeder}, \cite{kaletha1}, \cite{reeder}, the strategy of constructing a local Langlands correspondence is to first attach a character of a torus to a Langlands parameter, and then to construct a putative $L$-packet associated to this character.

We propose a different strategy to construct a local Langlands correspondence.  Langlands parameters and supercuspidal representations will be parameterized not by characters of tori, but by characters of covers of tori.  Given a Langlands parameter $\phi$, we first use a construction of Benedict Gross that naturally associates something close to a character $\chi$ of a cover of an elliptic torus.  We then attach a Weyl group orbit of conjectural supercuspidal characters to $\chi$ and prove that this orbit is the $L$-packet of supercuspidal representations that is associated to $\phi$ in \cite{debackerreeder}.  To do this, we rewrite supercuspidal characters in terms of functions on covers of elliptic tori as in Harish-Chandra's relative discrete series character formula.  We then compare our character formulas to those in \cite{debackerreeder}.

Let us briefly recall the construction of \cite{debackerreeder}.  To a certain class of Langlands parameters (tame, regular, semisimple, elliptic, Langlands parameters, or TRSELP for short) for an unramified connected reductive group $G$, DeBacker and Reeder associate a character of a torus $T$, to which they attach a collection of supercuspidal representations on the pure inner forms of $G(F)$, a conjectural $L$-packet.  To do this, DeBacker and Reeder extensively use the theory of Bruhat-Tits buildings.  They are also able to isolate the part of their $L$-packet corresponding to a particular pure inner form, and prove that their correspondence satisfies various natural conditions such as stability.

In the theory of real groups, an admissible homomorphism $W_{\mathbb{R}} \rightarrow {}^L G$ for a group $G(\mathbb{R})$ factors through the normalizer of a torus in ${}^L G$.  As such, it naturally produce a character $\chi$ of $T(\mathbb{R})_{\rho}$, where $T(\mathbb{R})_{\rho}$ is a double cover of some torus $T(\mathbb{R})$.  The cover $T(\mathbb{R})_{\rho}$ is called the $\rho$-\emph{cover} of $T(\mathbb{R})$.  Suppose $G(\mathbb{R})$ has relative discrete series representations.  Harish-Chandra has calculated the characters of the relative discrete series representations of $G(\mathbb{R})$.  Then the local Langlands correspondence for relative discrete series representations of $G(\mathbb{R})$ is given by attaching a Weyl group orbit of relative discrete series characters to $\chi$, and this is the motivation for our work over $p$-adic fields.

Very recently, Benedict Gross has shown that if one considers a Langlands parameter for a $p$-adic group that factors through the normalizer of a torus in ${}^L G$, then one obtains something close to a character of a cover of a torus, as follows.  Suppose that $G$ is a connected reductive group defined over a $p$-adic field $F$.  Let $\phi : W_F \rightarrow {}^L G$ be a Langlands parameter for $G(F)$, and suppose that $\phi$ factors through the normalizer of a torus. To $\phi$, one can associate an $F$-torus $T$ in $G$.  Suppose that $T$ splits over $E$, and set $\Gamma = Gal(E/F)$.  By the local Langlands correspondence for tori, we will show in section \ref{groupsoftypeL} that we can canonically associate to $\phi$ a character $\chi$ of $T(E)_{\Gamma}$, the group of coinvariants of $T(E)$ with respect to $\Gamma$.  Invariants and coinvariants are related by the norm map $$N : T(E) \rightarrow T(F)$$ $$t \mapsto \displaystyle\prod_{\xi \in \Gamma} \xi(t)$$ in the cohomology sequence
$$1 \rightarrow \hat{H}^{-1}(\Gamma,T(E)) \rightarrow T(E)_{\Gamma} \xrightarrow{N} T(F) = T(E)^{\Gamma} \rightarrow \hat{H}^0(\Gamma,T(E)) \rightarrow 1,$$ where $\hat{H}$ denotes Tate cohomology.

Suppose $\hat{H}^0(\Gamma, T(E))) = 0$, in which case $T(E)_{\Gamma}$ is then a cover of $T(F)$.  Let us also assume that $E/F$ is unramified.  We wish to attach a Weyl group orbit of conjectural supercuspidal characters to $\chi$, in analogy to the case of real groups.  These characters will be a $p$-adic analogue of the Harish-Chandra relative discrete series character formula.  Let $\Delta^+$ be a set of positive roots of $T(\overline{F})$ in $G(\overline{F})$.  Let $\rho$ denote half the sum of positive roots.  Let $\eta$ be any character of $E^*$ whose restriction to $F^*$ is a generator of $F^* / N_{E/F}(E^*)$, where $N_{E/F}$ is the norm map from $E$ to $F$.  Define $$\Delta^0(\gamma, \Delta^+) := \displaystyle\prod_{\alpha \in \Delta^+} (1 - \alpha^{-1}(\gamma)), \ \mathrm{for} \ \gamma \in T(F).$$

In the real case, the Weyl denominator of Harish-Chandra's relative discrete series character formula is given by $\Delta^0(\gamma, \Delta^+) \rho(\tilde{\gamma})$, where $\gamma \in  T(\mathbb{R})$ and $\tilde{\gamma}$ is any lift of $\gamma$ to $T(\mathbb{R})_{\rho}$ (by definition, $\rho$ is naturally a function on $T(\mathbb{R})_{\rho}$).  Since $\rho(\tilde{\gamma})^2 = 2 \rho(\gamma)$, $\rho(\tilde{\gamma})$ is canonically a square root of $2 \rho(\gamma)= 2 \rho(\Pi(\tilde{\gamma}))$, where $\Pi : T(\mathbb{R})_{\rho} \rightarrow T(\mathbb{R})$ is the canonical projection (see \S \ref{realgroupschapter}).  Our Weyl denominator will be a $p$-adic analogue of $\Delta^0(\gamma, \Delta^+) \rho(\tilde{\gamma})$, as we now explain.

If $\gamma \in T(F)$, note that $\Delta^0(\gamma, \Delta^+)$ does not take values in $\mathbb{C}^*$.  Since our character formula should take values in $\mathbb{C}^*$, we will apply $\eta$ to $\Delta^0(\gamma, \Delta^+)$.  For the $p$-adic version of $\rho(\tilde{\gamma})$, we wish to define a function on $T(E)_{\Gamma}$ that will act as a ``square root of $2 \rho(\gamma)$'', in analogy with the case of real groups.  Since our character formula should take values in $\mathbb{C}^*$, we want to define a function, denoted $\eta_{\rho}$, on $T(E)_{\Gamma}$, that acts as a ``square root of $\eta \circ 2 \rho$''.  More precisely, we want to define a function $\eta_{\rho}$ on $T(E)_{\Gamma}$ that is a canonical square root of $\eta \circ 2 \rho \circ N$.  Suppose that we have a way of defining $\eta_{\rho}$ in general.

We can now define our $p$-adic analogue of Harish-Chandra's relative discrete series character formula.  We first set $\overline{W} = N(G(E),T(E))^{\Gamma} / T(F)$ and $\overline{\overline{W}} = N(G(E),T(F)) / T(E)$, where $N(A,B)$ denotes the normalizer of $B$ in $A$.  It is not difficult to see that both $\overline{W}$ and $\overline{\overline{W}}$ act on $\widehat{T(E)_{\Gamma}}$.  We now define

\begin{equation}\label{characterformula}
\Theta_{\chi}(\gamma) = \epsilon(\chi, \Delta^+) \frac{ \displaystyle\sum_{n \in \overline{W}} \epsilon(n) n_* \chi(\tilde{\gamma}) }{\eta(\Delta^0(\gamma, \Delta^+)) \eta_{\rho}(\tilde{\gamma})}, \ \ \gamma \in T(F),
\end{equation}
where $n_* \chi(\tilde{\gamma}) := \chi(n^{-1} \tilde{\gamma} n)$ and where $\tilde{\gamma}$ is any element of $T(E)_{\Gamma}$ such that $N(\tilde{\gamma}) = \gamma$.  Here, $\epsilon(n) \in \{ \pm 1 \}$, and $\epsilon(\chi, \Delta^+)$ is some fixed constant, depending on $\chi$ and $\Delta^+$.  Note that we have assumed that $\hat{H}^0(\Gamma, T(E))) = 0$.  Therefore, $N$ is surjective, so there certainly exists $\tilde{\gamma} \in T(E)_{\Gamma}$ such that $N(\tilde{\gamma}) = \gamma$.

If $w \in \overline{\overline{W}}$, define ${}^w \Theta_{\chi}$ to be the same above formula, except that we replace $\chi$ by ${}^w \chi$, the $w$-conjugate of $\chi$.  Let $T(F)_{0,s}$ denote the set of strongly regular topologically semisimple elements of $T(F)$ (see \cite[\S 7]{debackerreeder}), and let $Z(F)$ denote the center of $G(F)$.

\begin{conjecture}\label{conjecture1}
\begin{enumerate}
\item ${}^w \Theta_{\chi}$ agrees with the character of a unique supercuspidal representation, denoted ${}^w \pi$, of $G(F)$ on $Z(F) T(F)_{0,s}$, where $Z(F)$ denotes the center of $G(F)$, and $T(F)_{0,s}$ denotes the set of strongly regular topologically semisimple elements of $T(F)$.
\item The assignment $$\phi \mapsto \{{}^w \pi : w \in \overline{\overline{W}} \}$$ is the local Langlands correspondence for (the pure inner forms of) $G(F)$.
\end{enumerate}
\end{conjecture}

We note that this conjecture holds for $GL(\ell,F)$, where $p > 2 \ell$ (see \cite{adrian1}).  In section \ref{characterformulas}, we will define a function $\eta_{\rho}$ on $T(E)_{\Gamma}$ that is a canonical square root of $\eta \circ 2 \rho \circ N$, in the case that $G(F) = PGSp(4,F)$. Our main theorem is the following.

\begin{theorem}\label{maintheorem}
The conjecture holds for depth zero supercuspidal $L$-packets of $PGSp(4,F)$.
\end{theorem}

We would like to note that in the case that $\hat{H}^0(\Gamma, T(E)) \neq 0$, the situation seems more difficult since our formula $\Theta_{\chi}$ is not defined on $Z(F) T(F)_{0,s}$ anymore, but namely on the image of $T(E)_{\Gamma}$ under the norm map.  However, one might be able to remedy this with a prediction of central character, as for example in \cite{grossreeder}.  We would also like to note that if $E/F$ is ramified, it is unclear how to define $\eta$, though this is the subject of future work.

We now present an outline of the paper.  In \S\ref{realgroupschater}, we provide the background material that we need from real groups as well as the motivation for our work.  In \S\ref{groupsoftypeL}, we review the theory of groups of type $L$ due to Gross, which is an integral part of our local Langlands correspondence.  In \S\ref{preliminaries}, we review the construction of \cite{debackerreeder}.  In \S\ref{rhoroots}, we conduct a thorough analysis of the tori that arise in the TRSELPs for $PGSp(4,F)$.  In \S\ref{cohomologyclasses}, we explicitly determine the structure of the genuine characters of $T(E)_{\Gamma}$.  In \S\ref{grossdebackerreeder}, we compare the construction of Gross to the construction of DeBacker/Reeder.  In \S\ref{characterformulas}, we define our conjectural character formula for depth zero supercuspidal representations of $PGSp(4,F)$.  We then prove theorem \ref{maintheorem}.

\section{Notation}\label{notation}

Suppose that $G$ is a connected reductive group over an arbitrary field $F$, and $T \subset G$ a torus defined over $F$. Let $\Delta^+$ be a set of positive roots of $G$ with respect to $T$.
We set $$\rho: = \frac{1}{2} \displaystyle\sum_{\alpha \in \Delta^+} \alpha.$$

Now let $F$ denote a nonarchimedean local field of characteristic zero.  We let $\mathfrak{o}_F$ denote the ring of integers of $F$, $\mathfrak{p}_F$ its maximal ideal, $\mathfrak{f}$ the residue field of $F$, $q$ the order of $\mathfrak{f}$, and $p$ the characteristic of $\mathfrak{f}$.  Let $\mathfrak{f_m}$ denote the degree $m$ extension of $\mathfrak{f}$.  We let $\varpi$ denote a uniformizer of $F$.  Let $F^u$ denote the maximal unramified extension of $F$.  Set $\Gamma_u = Gal(F^u/F)$.
We denote by $W_F$ the Weil group of $F$, $I_F$ the inertia subgroup of $W_F$, $I_F^+$ the wild inertia subgroup of $W_F$, and $W_F^{ab}$ the abelianization of $W_F$.  We denote by $W_F'$ the Weil-Deligne group, we set $W_t := W_F / I_F^+$, and we set $I_t := I_F / I_F^+$.  We fix an element $\Phi$ $\in Gal(\overline{F} / F)$ whose inverse induces the map $x \mapsto x^q$ on $\mathfrak{F} := \overline{\mathfrak{f}}$.

Now suppose $G$ is an unramified connected reductive group over $F$.  We fix $T \subset G$, an $F^u$-split maximal torus which is defined over $F$ and maximally $F$-split.  If $A$ and $B$ are groups, and if $B$ is a subgroup of $A$, then we let $N(A,B)$ denote the normalizer of $B$ in $A$.  We write $X := X_*(T)$, $W_o$ for the finite Weyl group $N(G(F^u), T(F^u)) / T(F^u)$, and set $N := N(G(F^u), T(F^u))$. We set $W = N^{Gal(F^u/F)} / T(F)$.  We denote by $\theta$ the automorphism of $X$ and $X \rtimes W_o$ induced by $\Phi$, where $X \rtimes W_o$ is the extended affine Weyl group.  For any finite Galois extension $E/F$, let $N_{E/F}:E \rightarrow F$ denote the norm map.  If $E/ F$ is a quadratic extension, we will sometimes denote the nontrivial Galois automorphism of $E/F$ by $x \mapsto \bar{x}$.  We let $E_1$ denote the quadratic unramified extension of $F$ and we let $E_2$ denote the quartic unramified extension of $F$.  Set $\Gamma_i = Gal(E_i/F)$.  If $A$ is a group and $B$ is a normal subgroup of $A$, we denote the image of $a \in A$ in $A / B$ by $[a]$.

\section{Background from real groups}\label{realgroupschapter}
\label{realgroupschater}

In order to motivate the theory that we wish to develop for $p$-adic groups, we describe the corresponding theory over $\mathbb{R}$ upon which our theory is based.  More information can be found in \cite{adamsvogan}.

\subsection{Covers of Tori}

It will be important to describe a part of the local Langlands correspondence having to do with relative discrete series representations, since this is the motivation for our construction over $p$-adic fields.

\begin{definition}\label{rhocoverdefinition}
Let $G$ be a connected reductive group over $\mathbb{R}$, $T \subset G$ a torus over $\mathbb{R}$, and let $\Delta^+$ be a set of positive roots of $G$ with respect to $T$.
Then $2 \rho \in X^*(T)$.  We define the \emph{$\rho$-cover $T(\mathbb{R})_{\rho}$ of $T(\mathbb{R})$} as the fiber product (in the category of groups) of the homomorphisms $2 \rho : T(\mathbb{R}) \rightarrow \mathbb{C}^*$ and the squaring map
$\Upsilon : \mathbb{C}^* \rightarrow \mathbb{C}^*$ given by $z \mapsto z^2$.  Thus, $T(\mathbb{R})_{\rho} = \{(t, \lambda) \in T(\mathbb{R}) \times \mathbb{C}^* : 2 \rho(t) = \lambda^2 \}$.
\end{definition}

Although $\rho$ is not necessarily a character of $T(\mathbb{R})$, it can naturally be thought of as a character of $T(\mathbb{R})_{\rho}$.  Namely, in the commutative diagram
$$
\begin{CD}
T(\mathbb{R})_{\rho} @> \Pi' >> \mathbb{C}^*\\
@VV \Pi V @VV \Upsilon V\\
T(\mathbb{R}) @>2 \rho>> \mathbb{C}^*
\end{CD}
$$
defining the fiber product, we have $\Pi'(\tilde{t})^2 = 2 \rho(t)$, where $\Pi(\tilde{t}) = t$.  Therefore, $\Pi'$ is a character of $T(\mathbb{R})_{\rho}$ which is a canonical square root of $2 \rho \circ \Pi$.  Throughout the rest of the paper, we will write $\rho$ instead of $\Pi'$.

The Weyl group acts on $T(\mathbb{R})_{\rho}$ as follows: If $(t, \lambda) \in T(\mathbb{R})_{\rho}$
and $s \in W(G(\mathbb{R}), T(\mathbb{R}))$, then define
\begin{equation}
s(t, \lambda) := (st, (s^{-1} \rho - \rho)(t) \lambda) \label{eq:weylgroupreal}
\end{equation}

\begin{definition}
A character $\tilde\chi : T(\R)_\rho\rightarrow \mathbb{C}^*$ is \emph{genuine} if it does not factor through $\Pi$.
\end{definition}

\begin{definition}
A genuine character
$\tilde\chi$ of $T(\mathbb{R})_{\rho}$ is called \emph{regular} if ${}^s \tilde\chi \neq \tilde\chi \ \forall s \in W(G(\mathbb{R}),T(\mathbb{R}))$ where ${}^s \tilde\chi(t, \lambda) := \tilde\chi(s^{-1}(t, \lambda))$.
\end{definition}

\subsection{Relative discrete series Langlands paramaters and character formulas for real groups}\label{realgroups}

In this section we will briefly describe the local Langlands correspondence for relative discrete series representations of real groups.  Let $G$ be a connected reductive group over $\mathbb{R}$ that contains a relatively compact maximal torus.  It is known that this is equivalent to $G(\mathbb{R})$ having relative discrete series representations.

\begin{definition}
Let $t$ be an indeterminate and let $k$ denote the rank of $G$.  For $h \in G(\R)$, define the Weyl denominator $D_G(h)$ by
$$\det(t + 1 - \Ad(h)) = D_G(h)t^k +\cdots  (\mbox{terms  of  higher  order})$$
Then if $\Delta$ is the set of roots of $T$ in $G$, $$D_G(h) = \prod_{\alpha \in \Delta} (1 - \alpha(h)).$$
\end{definition}

\begin{definition}
Let $G$ be a connected reductive group over $\mathbb{R}$, $T \subset G$ a maximal torus over $\mathbb{R}$.  Let $\Delta^+$ be a set of positive roots of $G$ with respect to $T$.  Define $$\Delta^0(h, \Delta^+) := \prod_{\alpha \in \Delta^+} (1 - \alpha^{-1}(h)), \ h \in T(\mathbb{R})$$
Then if the cardinality of $\Delta^+$ is $n$, we have $$(-1)^n D_G(h) = \Delta^0(h, \Delta^+)^2 (2 \rho)(h).$$
\end{definition}

Defining $|\rho(h)| := |2 \rho(h)|^{\frac{1}{2}}$ (the positive square root), we get that
$$|D_G(h)|^{\frac{1}{2}} = |\Delta^0(h, \Delta^+)| |\rho(h)|.$$
Although $|\rho|$ is a character of $T(\R)$, recall that $\rho$ is, in general, a character only of $T(\mathbb{R})_{\rho}$.  If $\tilde h \in T(\mathbb{R})_{\rho}$ maps to $h \in T(\mathbb{R})$ via the canonical projection, then $$|D_G(h)|^{\frac{1}{2}} = |\Delta^0(h, \Delta^+)| |\rho(h)| = |\Delta^0(h, \Delta^+)| |\rho(\tilde h)|.$$
We now present the classification of relative discrete series representations of $G(\mathbb{R})$.

\begin{theorem}\label{Harish-Chandra}(Harish-Chandra)
Let G be a connected reductive group, defined over $\mathbb{R}$.  Suppose that G contains a real Cartan subgroup T that is relatively compact.
Let $\tilde\chi$ be a genuine character of $T(\mathbb{R})_{\rho}$ that is regular.  Let $W := W(G(\mathbb{R}),T(\mathbb{R})) = N(G(\mathbb{R}),T(\mathbb{R})) / T(\mathbb{R})$
be the relative Weyl group.  Let $\epsilon(s) := (-1)^{\ell(s)}$ where $\ell(s)$ is the length of the Weyl group element $s \in W$. Let $T(\mathbb{R})\reg$ denote the regular set of $T(\mathbb{R})$.  Then there exists a unique constant $\epsilon(\tilde\chi, \Delta^+) = \pm 1$, depending only on $\tilde\chi$ and $\Delta^+$, and a unique relative discrete series representation of $G(\mathbb{R})$, denoted $\pi(\tilde\chi)$, such that
$$\theta_{\pi(\tilde\chi)}(h) =  \frac{\epsilon(\tilde\chi, \Delta^+)}{\Delta^0(h, \Delta^+) \rho(\tilde h)}\sum_{s \in W} \epsilon(s) \tilde\chi({}^s \tilde h), \quad \mbox{for $h \in T(\mathbb{R})\reg$}$$ where $\tilde h \in T(\mathbb{R})_{\rho}$ is any element such that $\Pi(\tilde h) = h$.  Moreover, every relative discrete series character of $G(\mathbb{R})$ is of this form.
\end{theorem}

It is important to note that while the numerator and denominator of the character formula live on $T(\mathbb{R})_{\rho}$, the quotient factors to a function on $T(\mathbb{R})\reg$.

We conclude the section by describing the local Langlands correspondence for relative discrete series representations of $G(\mathbb{R})$, assuming that $G(\mathbb{R})$ has relative discrete series. Fix a positive set of roots $\Delta^+$ of $G$ with respect to $T$.  Let $W_{\mathbb{R}}$ be the \emph{Weil group} of $\mathbb{R}$.  Let $\phi : W_{\mathbb{R}} \rightarrow {}^L G$
be a relative discrete series Langlands parameter.  The theory in \cite{adamsvogan} canonically attaches a genuine character $\tilde\chi$ of $T(\mathbb{R})_{\rho}$ to $\phi$.  Then the local Langlands correspondence for relative discrete series representations of $G(\mathbb{R})$ is given by attaching a Weyl group orbit of relative discrete series characters (as in theorem \ref{Harish-Chandra}) to $\tilde\chi$.  The rest of the paper will be devoted to proving the analogous result for $PGSp(4,F)$, where $F$ is a local non-Archimedean field of characteristic zero.

\section{Groups of type L}\label{groupsoftypeL}
We now review the theory of ``groups of type L'' due to Benedict Gross.  Let $F$ be a field, $F^{\mathrm s}$ a separable closure, and $T$ a torus defined over $F$ that splits over an extension $E \subset F^s$.
Let $\Gamma = \Gal(E/F)$.  Let $X^*(T)$ be the character module of $T$ and $X_*(T)$ the cocharacter
module of $T$.  Define $\hat{T} = X^*(T) \otimes \mathbb{C}^*$.
The group $\Gamma$ acts on $\hat{T}$ via its action on $X^*(T)$.

\begin{definition}
A \emph{group of type L} is a group extension of $\Gamma$ by $\hat{T}$.
\end{definition}

Let $D$ be such a group.  Then we have an exact sequence
$$1 \rightarrow \hat{T} \rightarrow D \rightarrow \Gamma \rightarrow 1$$
We now describe how, given a Langlands parameter $$\phi : W_F \rightarrow D,$$ where $D$ is a group of type L, we can naturally attach a character of $T(E)_{\Gamma} := T(E) / I_{\Gamma}(T(E))$, where $I_{\Gamma}(T(E)) = \{(1 - \xi)t \ : t \in T(E), \xi \in \Gamma \}$.
Restricting $\phi$ to $W_E$ we get a homomorphism $$\phi|_{W_E} : W_E \rightarrow \hat{T}$$
By the Langlands correspondence for tori, this gives us a character $\chi : T(E) \rightarrow \mathbb{C}^*$.  Since $\phi|_{W_E}$ extends to $\phi$, one can see that
$$\chi(t^{\sigma}) = \chi(t)\ \mbox{for all $\sigma \in \Gamma$.}$$
Therefore, $\chi(t^{\sigma - 1}) = 1$ for all $\sigma \in \Gamma$.  Thus, $\chi$ is trivial on the augmentation ideal $I_{\Gamma}(T(E))$
and gives $$\chi : T(E)_\Gamma \rightarrow \mathbb{C}^*$$
Invariants and coinvariants are related by the norm map $$N : T(E) \rightarrow T(F)$$ $$t \mapsto \displaystyle\prod_{\xi \in \Gamma} \xi(t)$$ in the Tate cohomology sequence
$$1 \rightarrow \hat{H}^{-1}(\Gamma,T(E)) \rightarrow T(E)_{\Gamma} \xrightarrow{N} T(F) = T(E)^{\Gamma} \rightarrow \hat{H}^0(\Gamma,T(E)) \rightarrow 1$$ (note that the norm map $N$ factors to $T(E)_{\Gamma}$).
We have thus constructed a character $\chi$ of $T(E)_{\Gamma}$ from a Langlands parameter $\phi$. We note that $T(E)_{\Gamma}$ is a cover of $N(T(E)_{\Gamma})$, which is a subgroup of $T(F)$.  It is sometimes the case that $N$ is surjective, in which case $\chi$ is then a character of $T(E)_{\Gamma}$, which is a cover of $T(F)$.

\section{Review of the Construction of $L$-packets of DeBacker and Reeder}\label{preliminaries}

We now review some of the basic theory from \cite{debackerreeder}.  Let $G$ be an unramified connected reductive group over $F$.  Let $\hat G$ denote the complex dual group of $G$.
Fix a pinning $(\hat{T}, \hat{B}, \{x_{\alpha} \})$ for $\hat G$ once and for all.
The operator $\hat{\theta}$ dual to $\theta$ extends to an automorphism of $\hat{T}$.  There is a unique extension of $\hat{\theta}$ to an automorphism of $\hat{G}$, satisfying $\hat{\theta}(x_{\alpha}) = x_{\theta \cdot \alpha}$ (see \cite[\S 3.2]{debackerreeder}).  Following \cite{debackerreeder}, we may form the semidirect product ${}^L G := \ \langle\hat{\theta}\rangle \ltimes \hat{G}$.

\begin{definition}
A Langlands parameter $\phi : W_F' \rightarrow {}^L G$
is called a \emph{tame regular semisimple elliptic Langlands parameter} (abbreviated TRSELP) if
\begin{enumerate}
\item $\phi$ is trivial on $I_F^+$,
\item The centralizer of $\phi(I_F)$
in $\hat{G}$ is a torus.
\item $C_{\hat{G}}(\phi)^o = (\hat{Z}^{\hat{\theta}})^o$, where $\hat{Z}$ denotes the center of $\hat{G}$, and where $C_{\hat{G}}(\phi)$ denotes the centralizer of $\phi$ in $\hat{G}$.
\end{enumerate}
\end{definition}

Condition (2) forces $\phi$ to be trivial on $SL(2,\mathbb{C})$.  Let $\hat{N} = N_{\hat{G}}(\hat{T})$.  After conjugating by an appropriate element of $\hat{G}$, we may assume that $\phi(I_F) \subset \hat{T}$ and $\phi(\Phi) = \hat{\theta} f$,
for some $f \in \hat{N}$.  Let $\hat{w}$ be the image of $f$ in $\hat{W}_o = \hat{N} / \hat{T}$,
and let $w$ be the corresponding element of $W_o$
under the natural identification of $W_o$ with $\hat W_o$.

Let $\phi$ be a TRSELP with associated $w$ and let $\sigma$ be the automorphism $w \theta$ of $T(F^u)$.
Let $\hat{\sigma}$ be the automorphism of $\hat T$
dual to $\sigma$, and let $n$ be the order of $\sigma$.  We set $\hat{G}\ab := \hat{G} / \hat{G}'$, where $\hat{G}'$ denotes the derived group of $\hat{G}$. Let ${}^L T_{\sigma} := \langle \hat{\sigma} \rangle \ltimes \hat{T}$.  DeBacker and Reeder (see \cite[\S 4]{debackerreeder}) associate to $\phi$ a $\hat{T}$-conjugacy class of Langlands parameters
\begin{equation}
\phi_T : W_t \rightarrow {}^L T_{\sigma} \ \label{phiT}
\end{equation}
as follows.  Set $\phi_T := \phi$ on $I_F$, and $\phi_T(\Phi) := \hat{\sigma} \ltimes \tau$ where $\tau \in \hat{T}$ is any element whose class in $\hat{T} / (1 - \hat{\sigma}) \hat{T}$ corresponds to the image of $f$ in $\hat{G}\ab / (1 - \hat{\theta}) \hat{G}\ab$ under the bijection

\begin{equation}
\hat{T} / (1 - \hat{\sigma}) \hat{T} \stackrel{\sim}{\rightarrow} \hat{G}\ab / (1 - \hat{\theta}) \hat{G}\ab \ \label{bijectionfortau}
\end{equation}

In \cite[Chapter 4]{debackerreeder}, DeBacker and Reeder construct a canonical bijection between $\hat{T}$-conjugacy classes of admissible homomorphisms $\phi : W_t \rightarrow {}^L T_{\sigma}$ and depth-zero characers of $T(F^u)^{\Phi_{\sigma}}$, where we identify
$T(F^u) = X \otimes F^u$ and $\Phi_{\sigma}$ is the automorphism $\sigma \otimes \Phi^{-1}$.
We briefly summarize this construction. Let $\mathbb{T} := X \otimes \mathfrak{F}^*$; note that $\Phi_{\sigma} = \sigma \otimes \Phi^{-1}$ also acts on $\mathbb{T}$.
Given automorphisms $\alpha, \beta$ of abelian groups $A,B$, respectively, let $\Hom_{\alpha, \beta}(A,B)$ denote the set of homomorphisms $f : A \rightarrow B$ such that $f \circ \alpha = \beta \circ f$.  The norm map $N_{\sigma}(t) = t \Phi_{\sigma}(t) \Phi_{\sigma}^2(t) \cdots \Phi_{\sigma}^{n-1}(t)$ induces isomorphisms $$\Hom(\mathbb{T}^{\Phi_{\sigma}}, \mathbb{C}^*) \stackrel{\sim}{\rightarrow}  \Hom_{\Phi_{\sigma}, Id}(\mathbb{T}^{\Phi_{\sigma}^n}, \mathbb{C}^*) \stackrel{\sim}{\rightarrow} \Hom_{\Phi_{\sigma}, Id}(X \otimes \mathfrak{f}_n^*, \mathbb{C}^*).$$
Moreover, the map $s \mapsto \chi_s$ gives an isomorphism
$$\Hom_{\Phi, \hat{\sigma}}(\mathfrak{f}_n^*, \hat{T}) \stackrel{\sim}{\rightarrow} \Hom_{\Phi_{\sigma}, Id}(X \otimes \mathfrak{f}_n^*, \mathbb{C}^*),$$
where $\chi_s(\lambda \otimes a) := \lambda(s(a))$.  The canonical projection $I_t \rightarrow \mathfrak{f}_m^*$ induces an isomorphism as $\Phi$-modules $I_t / (1 - Ad (\Phi)^m)I_t \stackrel{\sim}{\rightarrow} \mathfrak{f}_m^*$.  Since $\hat{\sigma}$ has order $n$, we have $\Hom_{\Phi, \hat{\sigma}}(\mathfrak{f}_n^*, \hat{T}) \cong \Hom_{Ad(\Phi), \hat{\sigma}}(I_t, \hat{T})$.  Therefore, the map $s \mapsto \chi_s$ is a canonical bijection $$\Hom_{Ad(\Phi), \hat{\sigma}}(I_t, \hat{T}) \stackrel{\sim}{\rightarrow} \Hom(\mathbb{T}^{\Phi_{\sigma}}, \mathbb{C}^*).$$
Moreover, we have an isomorphism
\begin{eqnarray*}
{}^0 T(F^u)^{\Phi_{\sigma}} \times X^{\sigma} \stackrel{\sim}{\rightarrow} T(F^u)^{\Phi_{\sigma}}\\
(\gamma, \lambda) \mapsto \gamma \lambda(\varpi),
\end{eqnarray*}
where ${}^0 T(F^u)$ denotes the maximal bounded subgroup of $T(F^u)$.

Finally, note that $\hat{T} / (1 - \hat{\sigma}) \hat{T}$ can be identified with the character group of $X^{\sigma}$, via the map taking $\tau \in \hat{T} / (1 - \hat{\sigma}) \hat{T}$ to $\chi_{\tau} \in \Hom(X^{\sigma}, \mathbb{C}^*)$, where $ \chi_{\tau}(\lambda) := \lambda(\tau)$.  Therefore, we have a canonical bijection between $\hat{T}$-conjugacy classes of admissible homomorphisms $\phi : W_t \rightarrow {}^L T_{\sigma}$ and depth-zero characters

\begin{equation}
\chi_{\phi} := \chi_s \otimes \chi_{\tau} \ \ \label{chitau}
\end{equation}

\noindent of $T(F^u)^{\Phi_{\sigma}}$, where $s := \phi|_{I_t}$, $\phi(\Phi) = \hat{\sigma} \ltimes \tau$.

\section{Generalities on Tori in $GSp(4,F)$ and $PGSp(4,F)$}\label{rhoroots}
In this section, we describe the unramified elliptic tori in $GSp(4,F)$ and $PGSp(4,F)$, since these are the tori that arise in \cite{debackerreeder}.  The torus
$$\dot S = \{\mathrm{diag}(x_1,x_2,x_3,x_4) : x_1 x_4 = x_2 x_3\}$$
is a split maximal tours of $GSp(4,\overline{F})$.
Its Weyl group has order $8$ and can be identified with the group of permutations of the elements
$\{x_1, x_2, x_3, x_4 \}$ that fix the relation $x_1 x_4 = x_2 x_3$.  We choose a system of positive roots
\begin{equation}
\label{eq:roots}
\alpha(x) = x_1 / x_2,\quad \beta(x) = x_2 / x_3,\quad (\alpha + \beta)(x) = x_1 / x_3,\quad (2 \alpha + \beta)(x) = x_1 / x_4,
\end{equation}
where $x = \mathrm{diag}(x_1,x_2,x_3,x_4)$.  All calculations in this paper will be performed with this choice of positive roots.  However, our results do not depend on this choice (see remark \ref{choiceofpositiveroots}).  We note that if $\rho$ denotes half the sum of the positive roots, then $2 \rho(x) = x_1^3/x_3^2 x_4$, where again $x = \mathrm{diag}(x_1,x_2,x_3,x_4)$.

Let $A = E_1 \times E_1$ and let $\sigma$ be the generator of $\Gal(E_1/F)$.  By \cite[\S1]{morris}, there is an elliptic maximal $F$-torus $\dot T_1$ of $GSp(4)$ such that
$\dot T_1(F)\cong \{ a \in A : a \sigma(a) \in F^* \}$, where $\sigma$ acts componentwise and $F^*$ is embedded diagonally as the set of pairs $\{(z,z)  : z \in F^*\}$.
More explicitly,
$$\dot T_1(F)  \cong \{ (x,y) \in E_1 \times E_1 : N_{E_1/F}(x) = N_{E_1/F}(y) \}.$$
We note that $\dot T_1$ is conjugate to $\dot S$ via an element $g\in GSp(4,E_1)$.  Moreover,
if $\mathrm{diag}(a,b,c,d)\in \dot S(E_1)$, $g$ can be chosen to satisfy
$\sigma (g^{-1} \mathrm{diag} (a,b,c,d) g) = g^{-1}\mathrm{diag} (\bar d,\bar c,\bar b,\bar a)g$.
Thus $\dot T_1(F) = \dot T_1(E_1)^{\Gal(E_1/F)}$ is $E_1$-conjugate to
the group of matrices of the form
$\mathrm{diag}(x,y,\bar y,\bar x)$, where $x,y\in E_1^*$ satisfy $N_{E_1/F}(x) = N_{E_1/F}(y)$.

Summarizing, we can and will identify
$\dot T_1(E_1)$ with
the group $\{(a,b,c,d) \in (E_1^*)^4 : ad = bc \}$, where $\sigma\in\Gamma_1$ acts via
the formula $\sigma (a,b,c,d) = (\bar{d},\bar{c},\bar{b},\bar{a})$ (recall that $\Gamma_i = Gal(E_i/F)$).  Then
$\dot{T}_1(F) = \dot{T}_1(E_1)^{\Gamma_1} = \{ (a, c, \bar{c}, \bar{a}) : a \bar{a} = c \bar{c} \}$,
which can further be identified with $\{ (x,y) \in E_1^* \times E_1^* : N_{E_1/F}(x) = N_{E_1/F}(y) \}$.

Again by \cite[\S 1]{morris}, there is an elliptic maximal $F$-torus $\dot T_2$ of $GSp(4)$ such that
$\dot T_2(F)\cong \{ x \in E_2 : x \tau^2(x) \in F^* \}$, where $\tau$ is a generator of $\Gal(E_2/F)$.
We will identify these groups in the following way.
We note that  $\dot T_2$ is conjugate to $\dot S$ via an element $h\in GSp(4,E_2)$.
In addition, if $(a,b,c,d)\in \dot S(E_2)$, then $h$ can be chosen to satisfy
$$\sigma (h^{-1} \mathrm{diag} (a,b,c,d) h) = h^{-1}\mathrm{diag} (\tau(c), \tau(a), \tau(d), \tau(b))h.$$
Thus $\dot T_2(F) = \dot T_2(E_2)^{\Gal(E_2/F)}$ is $E_2$-conjugate to
the group of matrices of the form
$$\mathrm{diag}(a,\tau(a), \tau^3(a), \tau^2(a)),$$
where $a\in E_2^*$ satisfies $a \tau^2(a)\in F^*$, i.e.,
$N_{E_2/E_1}(a)\in F^*$.

In summary, we may identify $\dot T_2(E_2)$ with
the group $\{(a,b,c,d) \in (E_2^*)^4 : ad = bc \}$, where $\tau\in\Gamma_2$ acts via
the formula $\tau (a,b,c,d) = (\tau(c),\tau(a),\tau(d),\tau(b))$.  Then
\begin{eqnarray*}
\dot{T}_2(F) &=& \dot{T}_2(E_2)^{\Gamma_2}\\
&=& \{ (a, \tau(a), \tau^3(a), \tau^2(a)) : a\in E_2^*, a \tau^2(a) = \tau (a\tau^2(a)) \}\\
&=& \{ (a, \tau(a), \tau^3(a), \tau^2(a)) : a\in E_2^*, a \tau^2(a)\in F^*\}
\end{eqnarray*}
We will further identify this group with $N_{E_2/E_1}^{-1}(F^*) =
\{ a\in E_2^*  : a\tau^2(a)\in F^* \}$.

\subsection{The torus $ T_1$.}\label{quadratictorus}
In this section, we will compute the Tate cohomology groups of $T_1$, and then give a simple description of the Tate cohomology exact sequence for $T_1$ (see section \ref{groupsoftypeL}).

Let $T_1$ be the image of $\dot T_1$ in $PGSp(4)$.  Then we can and will identify
$T_1(E_1)$ with
$\{(a,b,c,d)E_1^* : a,b,c,d \in E_1^*, ad=bc \}$, where $E_1^*$ is embedded diagonally.
The action of $\sigma$ on $\dot T_1(E_1)$ descends to an action on $T_1(E_1)$ given by
$\sigma( (a,b,c,d)E_1^* ) = (\bar{d},\bar{c},\bar{b},\bar{a})E_1^*$.
In order to compute Galois invariants and coinvariants of $T_1(E_1)$, it is useful to further identify $T_1(E_1)$ with $E_1^*\times E_1^*$ as follows.
\begin{lemma}
\label{T_1-identification}
Identifying $T_1(E_1)$ with $\{(a,b,c,d)E_1^* : a,b,c,d \in E_1^*, ad=bc \}$ as above, there is a Galois-equivariant isomorphism
$$\phi : T_1(E_1)  \rightarrow  E_1^* \times E_1^*$$
$$(a,b,c,d)E_1^*  \mapsto  (a/b, b/c)$$
where the Galois action on $E_1^* \times E_1^*$ is given by $\sigma(w,z) \mapsto (1/\bar{w}, 1/\bar{z})$.
Thus
$$T_1(F) = T_1(E_1)^{\Gamma_1}\cong \ker(N_{E_1/F})\times \ker(N_{E_1/F}).$$
\end{lemma}

\proof
The map is easily seen to be well defined.  To show injectivity, suppose $\phi((a,b,c,d)E_1^*) = (1,1)$.  Thus, $a = b = c$.  Since $ad = bc$, we have $a = b = c = d$, so $(a,b,c,d)E_1^*$ is trivial in $T_1(E_1)$.  The inverse of $\phi$ is given by $(w,z) \mapsto (w,1,1/z,1/wz)E_1^*$.  To show that $\phi$ is Galois equivariant, note that
\begin{eqnarray*}
\phi(\sigma((a,b,c,d)E_1^*)) &=& \phi( (\bar{d},\bar{c},\bar{b},\bar{a})E_1^*)\\
&=& (\bar{d}/\bar{c}, \bar{c}/\bar{b})\\
&=& (\bar{b}/\bar{a}, \bar{c}/\bar{b})\\
&=& \sigma( a/b,b/c) \\
&=& \sigma(\phi( (a,b,c,d)E_1^*))
\end{eqnarray*}
since $ad=bc$.  The final statement now follows easily.
\qed

We now consider the norm map $N : T_1(E_1) \rightarrow T_1(E_1)^{\Gamma_1}$ given in section \ref{groupsoftypeL}.
\begin{lemma}\label{T_1-H0}
$\hat{H}^0(\Gamma_1, T_1(E_1)) = 0$.
\end{lemma}

\proof
Recall that $\hat{H}^0(\Gamma_1, T_1(E_1)) = T_1(E_1)^{\Gamma_1} / N(T_1(E_1))$, where $N$ denotes the norm map to $T_1(F)$.
Identifying $T_1(E_1)$ with $E_1^*\times E_1^*$ as in Lemma~\ref{T_1-identification}, we have $N(w,z) = (w/\bar w,z/\bar z)$.  Thus, by Hilbert's Theorem 90, $N(T_1(E_1)) = \ker(N_{E_1/F})\times \ker(N_{E_1/F}) = T_1(E_1)^{\Gamma_1}$.
\qed

\begin{lemma}\label{T_1-H-1}
$\hat{H}^{-1}(\Gamma_1, T_1(E_1)) \cong (\mathbb{Z} / 2 \mathbb{Z})\times (\mathbb{Z} / 2 \mathbb{Z})$.
\end{lemma}

\proof
Let $1-\sigma: T_1(E_1)\longrightarrow T_1(E_1)$ be the map $x\mapsto x/\sigma(x)$.  We note that $\hat{H}^{-1}(\Gamma_1, T_1(E_1)) = \ker(N)/ \im (1-\sigma)$.  Let $(w,z)\in E_1^*\times E_1^*$, which we identify with $T_1(E_1)$ as above.  Since $N(w,z) = (w/\bar w,z/\bar z)$, it follows that $\ker (N) = F^*\times F^*$.

On the other hand, for $(x,y)\in E_1^*\times E_1^* = T_1(E_1)$, we have
$(1-\sigma)(x,y) = (N_{E_1/F}(x),N_{E_1/F}(y))$,
so $\im(1-\sigma) = N_{E_1/F}(E_1^*)\times N_{E_1/F}(E_1^*)$.
Thus
$$\ker(N)/\im(1-\sigma) = (F^*/(N_{E_1/F}(E_1^*))\times (F^*/N_{E_1/F}(E_1^*))\cong (\mathbb{Z} / 2 \mathbb{Z})\times (\mathbb{Z} / 2 \mathbb{Z}).$$
\qed

The map $N$ factors through the group $T_1(E_1)_{\Gamma_1}$ of coinvariants to give a map $N : T_1(E_1)_{\Gamma_1} \rightarrow T_1(E_1)^{\Gamma_1}$.  We now give a concrete description of $T_1(E_1)_{\Gamma_1}$.
For $(w,z) \in E_1^* \times E_1^*$, $(1-\sigma)(w,z) = (w,z)(1/\bar{w},1/\bar{z})^{-1} = (N_{E_1/F}(w), N_{E_1/F}(z))$.  Therefore, $T_1(E_1)_{\Gamma_1}$ can be identified with
$$\frac{E_1^* \times E_1^*}{(1-\sigma)(E_1^* \times E_1^*)} = \left(E_1^* / N_{E_1/F}(E_1^*)\right) \times \left(E_1^* / N_{E_1/F}(E_1^*)\right).$$
Recalling that $T_1(E_1)^{\Gamma_1} = \ker(N_{E_1/F})\times \ker(N_{E_1/F})$, we have that
$N: T_1(E_1)_{\Gamma_1} \rightarrow T_1(E_1)^{\Gamma_1}$ is given by
$$(wN_{E_1/F}(E_1^*),wN_{E_1/F}(E_1^*))\mapsto (w/\bar w,z/\bar z).$$
Therefore, the exact sequence of Tate cohomology groups in section \ref{groupsoftypeL} reduces in our case to the standard exact sequence
$$1 \rightarrow \mathbb{Z} / 2 \mathbb{Z} \times \mathbb{Z} / 2 \mathbb{Z} \rightarrow E_1^* / N_{E_1/F}(E_1^*) \times E_1^* / N_{E_1/F}(E_1^*) \xrightarrow{N} \ker(N_{E_1/F})\times \ker(N_{E_1/F}) \rightarrow 1$$

\subsection{The torus $T_2$}\label{quartictorus}
In this section, we will compute the Tate cohomology groups of $T_2$, and then give a simple description of the Tate cohomology exact sequence for $T_2$ (see section \ref{groupsoftypeL}).

Let $T_2$ be the image of $\dot T_2$ in $PGSp(4)$.  Then we can and will identify
$T_2(E_2)$ with
$\{(a,b,c,d)\in E_2^* : a,b,c,d \in E_2^*, ad=bc \}$, where $E_2^*$ is embedded in diagonally.
The action of $\tau$ on $\dot T_2(E_2)$ descends to an action on $T_2(E_2)$ given by
$\tau( (a,b,c,d)E_2^* ) = (\tau(c),\tau(a),\tau(d),\tau(b))E_2^*$.
In order to compute Galois invariants and coinvariants of $T_2(E_2)$, the following further identification is useful.

\begin{lemma}
\label{T2-identification}
Identifying $T_2(E_2)$ with $\{(a,b,c,d)E_2^* : a,b,c,d \in E_2^*, ad=bc \}$ as above, there is a Galois-equivariant isomorphism
\begin{eqnarray*}
\varphi : T_2(E_2) \xrightarrow{\sim} E_2^* \times E_2^*\\
(a,b,c,d)E_2^* \mapsto (a/b, a/c),
\end{eqnarray*}
where the Galois action on $E_2^* \times E_2^*$ is given by $\tau(x,y) := (\tau(y)^{-1}, \tau(x))$.
\end{lemma}

\proof
The map is clearly well defined, injective, and surjective as in the proof of Lemma~\ref{T_1-identification}.  Moreover, if $(a,b,c,d)E_2^* \in T_2(E_2)$ (so that $ad = bc$), then
\begin{eqnarray*}
\varphi(\tau((a,b,c,d)E_2^*)) &=& \varphi((\tau(c),\tau(a),\tau(d),\tau(b))E_2^*)\\
&=& (\tau(c/a),\tau(c/d))\\
&=& (\tau(c/a),\tau(a/b))\\
&=& \tau(a/b,a/c)\\
&=& \tau(\varphi((a,b,c,d)E_2^*)).
\end{eqnarray*}
\qed



\begin{lemma}\label{T_2-H0}
$\hat{H}^0(\Gamma_2, T_2(E_2)) = 0$
\end{lemma}

\proof
Recall that $\hat{H}^0(\Gamma_2, T_2(E_2)) = T_2(E_2)^{\Gamma_2} / N(T_2(E_2))$, where $N$ denotes the norm map to $T_2(F)$ and $\Gamma_2 = \Gal (E_2/F)$.  In showing that this quotient is trivial, we will identify $T_2(E_2)$ with $E_2^*\times E_2^*$ as in Proposition~\ref{T2-identification}.
Suppose that $(x,y) \in T_2(E_2)^{\Gamma_2}$ (under the above identification).  Then one easily sees that $y = \tau(x)$ and $x\in\ker(N_{E_2/E_1})$.  Since $\tau^2$ generates $\Gal (E_2/E_1)$, it follows from Hilbert's Theorem 90 that $x = w/\tau^2(w)$ for some $w\in E_2^*$.  Then
\begin{eqnarray*}
N(w,1) &=& (w,1)\cdot \tau(w,1) \cdot \tau^2(w,1)\cdot \tau^3(w,1)\\
&=& (w,1)\cdot (1,\tau(w))\cdot (\tau^2(w)^{-1},1)\cdot (1,\tau^3(w)^{-1})\\
&=& (w/\tau^2(w),\tau(w/\tau^2(w)))\\
&=& (x,\tau (x))
\end{eqnarray*}
It follows that $(x,y)\in N(T_2(E_2))$, so $N(T_2(E_2)) = T_2(E_2)^{\Gamma_2}$.
\qed



\begin{lemma}\label{T_2-H-1}
$\hat{H}^{-1}(\Gamma_2, T_2(E_2)) \cong \mathbb{Z} / 2 \mathbb{Z}$
\end{lemma}

\proof
Let $\vartheta: T_2(E_2)\longrightarrow T_2(E_2)$ be the map $x\mapsto x/\tau(x)$.  We note that $\hat{H}^{-1}(\Gamma_2, T_2(E_2)) = \ker(N)/ \im (\vartheta)$.  Let $(w,z)\in E_2^*\times E_2^*$, which we identify with $T_2(E_2)$ as above.  Note that
\begin{eqnarray*}
N(w,z) &=& (w,z)\cdot \tau(w,z) \cdot \tau^2(w,z)\cdot \tau^3(w,z)\\
&=& (w,z)\cdot (\tau(z)^{-1},\tau(w))\cdot (\tau^2(w)^{-1},\tau^2(z)^{-1})\cdot (\tau^3(z),\tau^3(w)^{-1})\\
&=& (w\tau(z)^{-1}\tau^2(w)^{-1}\tau^3(z),z\tau(w)\tau^2(z)^{-1}\tau^3(w)^{-1}) .
\end{eqnarray*}
Thus $N(w,z) = (x,\tau(x))$, where $x = w\tau(z)^{-1}\tau^2(w)^{-1}\tau^3(z)$.  Thus the condition $N(w,z) = (1,1)$
is equivalent to $x=1$, or $w/\tau(z) = \tau^2(w)/\tau^3(z) = \tau^2(w/\tau(z))$.  It follows that $w/\tau(z)\in  E_1^*$, so that $\tau (z) = w\alpha$ for some $\alpha\in E_1^*$.  Thus
$$\ker (N) = \{ (w,\tau^{-1}(w)\beta : w\in E_2^*, \beta\in E_1^*\} .$$

For $(x,y)\in E_2^*\times E_2^* = T_2(E_2)$, we have
$$\vartheta(x,y) = (x,y)\tau(x,y)^{-1} = (x,y)(\tau(y),\tau(x)^{-1}) = (x\tau(y),y/\tau(x)).$$
Making a change of variables by setting $w = x \tau(y)$, we get that
$$y / \tau(x) = y \tau^2(y) / \tau(w) = \tau^{-1}(w)\cdot\frac{N_{E_2/E_1}(y)}{\tau(w)\tau^{-1}(w)} =
\tau^{-1}(w)\cdot N_{E_2/E_1}(y/\tau(w)).$$
It follows that
\begin{equation}
\label{eq:augmentation}
\im(\vartheta) = \{ (w, \tau^{-1}(w)\beta : w\in E_2^*, \beta\in N_{E_2/E_1}(E_2^*) \}.
\end{equation}
Therefore, it is clear that $\ker(N)/\im(\vartheta) \cong \Z/2\Z$.
\qed

\begin{lemma}
\label{other}
There are natural isomorphisms
\begin{eqnarray*}
T_2(E_2)_{\Gamma_2} &\cong& E_2^* / N_{E_2/E_1}(E_2^*)\\
T_2(E_2)^{\Gamma_2}  &\cong& \ker(N_{E_2/E_1}).
\end{eqnarray*}
Identifying these pairs of isomorphic groups, the norm map
$N: T_2(E_2)_{\Gamma_2} \rightarrow T_2(E_2)^{\Gamma_2}$ is given by
$$wN_{E_2/E_1}(E_2^*)\mapsto w/\tau^2(w)\quad\mbox{for $w\in E_2^*$}.$$
\end{lemma}

\proof
Identify $T_2(E_2)$ with $E_2^*\times E_2^*$ as above.  Using (\ref{eq:augmentation}), we see that the surjective homomorphism $T_2(E_2)\rightarrow E_2^*$ defined by $(w,z)\mapsto w \tau(z)^{-1}$ maps the augmentation
ideal $\im (\vartheta)$ onto $N_{E_2/E_1}(E_2^*)$.  It therefore defines an isomorphism $T_2(E_2)_{\Gamma_2}\cong E_2^* / N_{E_2/E_1}(E_2^*)$.  Moreover, by the proof of Lemma~\ref{T_2-H0}, we have
$$T_2(E_2)^{\Gamma_2} = \{(x,\tau(x)):x\in\ker(N_{E_2/E_1})\} \cong \ker(N_{E_2/E_1})$$
via the map $(x,\tau(x))\mapsto x$.

Let $(w,1)\in E_2^* \times E_2^*$, which we identify with $T_2(E_2)$ as above.  Let $\overline{(w,1)}$ be the image
of $(w,1)$ in $T_2(E_2)_{\Gamma_2}$.  As in the proof of Lemma \ref{T_2-H0},
$$N\overline{(w,1)} = (w/\tau^2(w),\tau(w/\tau^2(w))).$$
Under the isomorphisms in the preceding paragraph, $\overline{(w,1)}\in T_2(E_2)_{\Gamma_2}$ corresponds to
$wN_{E_2/E_1}(E_2^*)\in E_2^*/N_{E_2/E_1}(E_2^*)$, while $(w/\tau^2(w),\tau(w/\tau^2(w)))\in T_2(E_2)^{\Gamma_2}$ corresponds to $w/\tau^2(w)\in \ker(N_{E_2/E_1})$.
\qed

Therefore, we have finally shown that the exact sequence of Tate cohomology groups in section \ref{groupsoftypeL} reduces therefore in our case to the standard exact sequence
$$1 \rightarrow \mathbb{Z} / 2 \mathbb{Z} \rightarrow E_2^* / N_{E_2/E_1}(E_2^*) \xrightarrow{N} \ker(N_{E_2/E_1}) \rightarrow 1$$

\section{Cohomology class of group of type L}\label{cohomologyclasses}
Denote by $\sigma$ the generator of $Gal(E_1/F)$ and by $\tau$ a generator of $Gal(E_2/F)$.
Let $G = PGSp(4)$.  Let $\phi : W_F \rightarrow {}^L G$ be a TRSELP for $G(F)$.  Note that ${}^L G = \Spin(5, \mathbb{C}) \times \Gal(\overline{F}/F)$.

Recall that by the theory of groups of type L, $\phi$ naturally gives rise to a character $\chi$ of the group of coinvariants $T(E)_{\Gamma}$ of a torus $T$ that is defined over $F$ and split over $E$, where $\Gamma = Gal(E/F)$.  In this section, we explicitly compute the restriction of $\chi$ to $\hat{H}^{-1}(\Gamma, T(E))$ (see section \ref{groupsoftypeL}).

Given a root $\alpha_*$ of $\hat T$ in ${\rm Spin}(5)$, we denote by $w_{\alpha_*}$ the reflection in the Weyl group of $\hat T$ in ${\rm Spin}(5)$ through $\alpha_*$.
For each such $\alpha_*$, there is an associated homomorphism
${\rm SL}(2) \rightarrow{\rm Spin}(5)$.  Let $n_{\alpha_*}$ be the image of the matrix
$\left(
\begin{smallmatrix}
0 & 1\\
-1 & 0
\end{smallmatrix}
\right) $
under this map.  Then $n_{\alpha_*}$ lies in the normalizer $N(\hat{G}, \hat{T})$ of $\hat T$ in ${\rm Spin}(5)$, and its image in the Weyl group of $\hat T$ is $w_{\alpha_*}$.  Let $\alpha_*^\vee$ denote the coroot of $\hat T$ corresponding to $\alpha_*$.  Then
\begin{equation}
\label{eq:square}
n_{\alpha_*}^2 = \alpha_*^\vee(-1).
\end{equation}
We have the following commutation relations for all roots $\alpha_*,\beta_*$ and $t\in\C$ (see~\cite[\S7.2]{carter})
\begin{eqnarray}
n_{\alpha_*} \beta_*^\vee(t) n_{\alpha_*}^{-1} &=& (w_{\alpha_*}(\beta_*))^\vee(t)\label{eq:chevalley1}\\
n_{\alpha_*} n_{\beta_*} n_{\alpha_*}^{-1} &=& (w_{\alpha_*}(\beta_*))^\vee(\eta_{\alpha_*,\beta_*})\cdot n_{w_{\alpha_*}(\beta_*)}\quad
\mbox{for a certain $\eta_{\alpha_*,\beta_*}\in\C$}\label{eq:chevalley2}.
\end{eqnarray}

We will henceforth denote by $\alpha_*$ and $\beta_*$, respectively, the long and short simple roots of $\hat T$ in ${\rm Spin}(5)$ (with respect to a fixed pinning).

\begin{lemma}
\label{tits-group}
Let $\hat n = n_{\alpha_*} n_{\beta_*} n_{\alpha_*} n_{\beta_*}$.  Then $\hat n^2 = \beta_*^\vee(-1)$.
\end{lemma}
\proof
Observe that
\[
\hat n^2 = n_{\alpha_*} n_{\beta_*} n_{\alpha_*} n_{\beta_*} n_{\alpha_*} n_{\beta_*} n_{\alpha_*} n_{\beta_*}
\]
\[
= n_{\alpha_*} n_{\beta_*} \cdot (n_{\alpha_*} n_{\beta_*} n_{\alpha_*}^{-1}) \cdot n_{\alpha_*}^2 \cdot n_{\beta_*} n_{\alpha_*} n_{\beta_*}
\]
\[
= n_{\alpha_*} n_{\beta_*}\cdot  (\alpha_*+\beta_*)^\vee(\eta_{\alpha_*,\beta_*})\cdot n_{\alpha_*+\beta_*}\cdot  \alpha_*^\vee(-1) \cdot n_{\beta_*} n_{\alpha_*} n_{\beta_*}
\]
by (\ref{eq:square}) and (\ref{eq:chevalley2}), since $w_{\alpha_*}(\beta_*) = \alpha_*+\beta_*$.
The above expression is equal to
\begin{eqnarray}
\label{eq:coroots}
\lefteqn{n_{\alpha_*}\cdot  \left(n_{\beta_*}\left((\alpha_*+\beta_*)^\vee(\eta_{\alpha_*,\beta_*})\cdot n_{\alpha_*+\beta_*}\cdot  \alpha_*^\vee(-1) \right)n_{\beta_*}^{-1}\right) \cdot n_{\beta_*}^2 \cdot n_{\alpha_*} n_{\beta_*}}\hspace{2cm}\nonumber\\
&=& n_{\alpha_*}\cdot  \left((\alpha_*+\beta_*)^\vee(\eta_{\alpha_*,\beta_*})\cdot (\alpha_*+\beta_*)^\vee(\eta_{\beta_*,\alpha_*+\beta_*})
\cdot n_{\alpha_*+\beta_*}\cdot  (\alpha_*+2\beta_*)^\vee(-1) \right)\nonumber \\
& & \cdot \beta_*^\vee (-1) \cdot n_{\alpha_*} n_{\beta_*}
\end{eqnarray}
again by (\ref{eq:square}), (\ref{eq:chevalley1}), and (\ref{eq:chevalley2}), since $w_{\beta_*}(\alpha_*+\beta_*) = \alpha_*+\beta_*$ and
$w_{\beta_*}(\alpha_*) = \alpha_*+2\beta_*$.
But (\ref{eq:coroots}) can be rewritten as
\begin{eqnarray}
\label{eq:coroots2}
\lefteqn{n_{\alpha_*} \left((\alpha_*+\beta_*)^\vee(\eta_{\alpha_*,\beta_*}\ \eta_{\beta_*,\alpha_*+\beta_*})
\cdot n_{\alpha_*+\beta_*}\cdot  (\alpha_*+2\beta_*)^\vee(-1) \cdot \beta_*^\vee (-1)\right)n_{\alpha_*}^{-1} \cdot n_{\alpha_*}^2 \cdot n_{\beta_*}}\hspace{2cm}\nonumber\\
&=& \beta_*^\vee(\eta_{\alpha_*,\beta_*}\ \eta_{\beta_*,\alpha_*+\beta_*})
\cdot \beta_*^\vee(\eta_{\alpha_*,\alpha_*+\beta_*})\cdot n_{\beta_*}\cdot  (\alpha_*+2\beta_*)^\vee(-1)\nonumber \\
& & \cdot (\alpha_*+\beta_*)^\vee (-1) \cdot \alpha_*^\vee (-1)\cdot n_{\beta_*} ,
\end{eqnarray}
since
$$w_{\alpha_*}(\alpha_* + \beta_*) = \beta_*, \quad w_{\alpha_*}(\alpha_* + 2\beta_*) = \alpha_*+2\beta_*, \quad
w_{\alpha_*}(\beta_*) = \alpha_*+\beta_* .$$
We may re-express (\ref{eq:coroots2}) as
\begin{eqnarray}
\lefteqn{\beta_*^\vee(\eta_{\alpha_*,\beta_*}\ \eta_{\beta_*,\alpha_*+\beta_*}\ \eta_{\alpha_*,\alpha_*+\beta_*})\cdot n_{\beta_*}\left(  (\alpha_*+2\beta_*)^\vee(-1) \cdot (\alpha_*+\beta_*)^\vee (-1) \cdot \alpha_*^\vee (-1)\right)n_{\beta_*}^{-1}\cdot n_{\beta_*}^2}\nonumber\\
&=& \beta_*^\vee(\eta_{\alpha_*,\beta_*}\ \eta_{\beta_*,\alpha_*+\beta_*}\ \eta_{\alpha_*,\alpha_*+\beta_*})\cdot \alpha_*^\vee(-1) \cdot (\alpha_*+\beta_*)^\vee (-1) \cdot (\alpha_*+2\beta_*)^\vee (-1)\cdot \beta_*^\vee(-1),\nonumber
\end{eqnarray}
as above, noting that $w_{\beta_*}(\alpha_*+2\beta_*) = \alpha_*$.  One checks easily that $(\alpha_*+\beta_*)^\vee = 2\alpha_*^\vee + \beta_*^\vee$ and $(\alpha_*+2\beta_*)^\vee = \alpha_*^\vee + \beta_*^\vee$.  Therefore, the preceding displayed expression equals
\begin{eqnarray}
\lefteqn{\beta_*^\vee(-\eta_{\alpha_*,\beta_*}\ \eta_{\beta_*,\alpha_*+\beta_*}\ \eta_{\alpha_*,\alpha_*+\beta_*})\cdot \alpha_*^\vee(-1) \cdot \left[ 2\alpha_*^\vee (-1)\cdot\beta_*^\vee(-1)\right] \cdot \left[\alpha_*^\vee (-1)\beta_*^\vee(-1)\right]}\hspace{9cm}\nonumber\\
&=& \beta_*^\vee(-\eta_{\alpha_*,\beta_*}\ \eta_{\beta_*,\alpha_*+\beta_*}\ \eta_{\alpha_*,\alpha_*+\beta_*})\nonumber
\end{eqnarray}
Basic results on Chevalley groups (see \cite{carter}) allow one to compute that
\[
-\eta_{\alpha_*,\beta_*}\ \eta_{\beta_*,\alpha_*+\beta_*}\ \eta_{\alpha_*,\alpha_*+\beta_*} = -1,
\]
concluding the proof.
\qed

Recall that our parameter satisfies $\phi(I_F) \subset \hat{T}$ and $\phi(\Phi) = \hat{\theta} f$, for some $f \in N_{\hat{G}}(\hat{T})$.  The Weyl group of $Spin(5)$ is the dihedral group of order $8$.  If we set $r = w_{\alpha_*} w_{\beta_*}$, then for $\phi$ to be a TRSELP, it is necessary that $\hat{w}$ must be a nontrivial element of $\langle r\rangle$.
Therefore, without loss of generality, we may take $\hat{w}$ to be $r$ or $r^2$ (since $r^3$ is conjugate to $r$ in the Weyl group of $Spin(5)$, so replacing $r$ by $r^3$ would yield an equivalent TRSELP).

\begin{definition}\label{trselptype}
We say that a TRSELP $\phi$ is \emph{of type $r$} if $\hat{w} = r$.  We say that a TRSELP $\phi$ is \emph{of type $r^2$} if $\hat{w} = r^2$.
\end{definition}

Let $\phi$ be a TRSELP of type $r^2$.  We first note that this type of TRSELP is associated to the torus $T_1$.
More precisely, the Weyl group element associated to $\phi$ acts on the module $X^*(T_1) = X_*(\hat{T_1})$ by $(a,b) \mapsto (-a,-b)$, which is exactly the action that gives rise to $T_1$.   We have identified $T_1(E_1)$ with $E_1^* \times E_1^*$ in section \ref{quadratictorus} via the roots $\alpha$ and $\beta$ of $PGSp(4)$.  More precisely, our identification was

$$\varphi : T_1(E_1) \xrightarrow{\sim} E_1^* \times E_1^*$$
$$(a,b,c,d)E_1^* \mapsto (a/b, b/c) = (\alpha((a,b,c,d)E_1^*), \beta((a,b,c,d)E_1^*)).$$
This gives us an identification of $\hat{T}_1$ with $\mathbb{C}^* \times \mathbb{C}^*$ that is compatible with duality of tori.  Namely, we use the images of the simple roots $\alpha$ and $\beta$ of $PGSp(4)$ under the isomorphism $X^*(T_1) \xrightarrow{\sim} X_*(\hat{T}_1)$ to decompose $\hat{T}_1$.  Explicitly, the map $\mathbb{C}^* \times \mathbb{C}^* \xrightarrow{\sim} \hat{T}_1$ is given by
\[
(w,z) \mapsto \alpha_*^{\vee}(w)\beta_*^{\vee}(z).
\]
Moreover, we have that $\Gal (E_1/F)$ acts on $\hat{T}_1 = \C^*\times \C^*$ via $\sigma(w,z) = (1/w, 1/z)$.  By the theory in section \ref{groupsoftypeL}, $\phi$ naturally gives rise to a character $\chi$ of $T_1(E_1)_{\Gamma_1}$.  We wish to explicitly compute $\chi|_{\hat{H}^{-1}(\Gamma_1, T_1(E_1))}$.

Recall the exact sequence $$1 \rightarrow \hat{H}^{-1}(\Gamma_1, T_1(E_1)) \rightarrow T_1(E_1)_{\Gamma_1} \rightarrow T_1(F) \rightarrow 1$$
(since $\hat{H}^0(\Gamma_1, T_1(E_1)) = 1$).  Recall that $\phi|_{W_{E_1}}$ has image in $\hat{T}_1$.  Therefore, since the abelianization of $W_{E_1}$ is $E_1^*$, $\phi|_{W_{E_1}}$ factors through a map $\phi' : E_1^* \rightarrow \hat{T}_1$.  Recall our identification of $\hat{T}_1$ with $\mathbb{C}^* \times \mathbb{C}^*$ via $\alpha_*^{\vee}$ and $\beta_*^{\vee}$.  We therefore get a canonical map
$$E_1^* \rightarrow \mathbb{C}^* \times \mathbb{C}^*$$
$$x \mapsto (\chi_1(x), \chi_2(x))$$
The local Langlands correspondence for tori says that $\chi(w,z) = \chi_1(w) \chi_2(z)$ for all $w, z \in E_1^*$, where we view $\chi$ as a character of $E_1^* \times E_1^*$ via the identification of $T_1(E_1)$ with $E_1^* \times E_1^*$ in section \ref{quadratictorus}. Under this identification of $T_1(E_1)$ with $E_1^* \times E_1^*$, $\hat{H}^{-1}(\Gamma_1, T_1(E_1))$  is identified canonically with $F^* / N_{E_1/F}(E_1^*) \times F^* / N_{E_1/F}(E_1^*)$.  Therefore, we wish to compute $\chi(\varpi, \varpi), \chi(\varpi, 1)$, and $\chi(1, \varpi)$.

Recall that the Artin map $W_{E_1} \rightarrow E_1^*$ sends $\Phi^2$ to $\varpi$. Write $\phi(\Phi^2) = \phi'(\varpi)$ as $\alpha_*^{\vee}(a) \beta_*^{\vee}(b)$, for some $a,b \in \mathbb{C}^*$.  The above discussion allows us to conclude that $\chi(\varpi, \varpi) = ab$, $\chi(\varpi,1) = a$, and $\chi(1, \varpi) = b$.

\begin{proposition}
\label{prop:trselp1}
We have $\chi(\varpi, \varpi) = -1, \chi(\varpi,1) = 1$, and $\chi(1,\varpi) = -1$.
Therefore, $\chi$ factors to a genuine character of the two-fold cover
$\widetilde{T_1(F)} := T_1(E_1)_{\Gamma_1} / \langle (\varpi,1) \rangle$ of $T_1(F)$.
\end{proposition}

\proof

By the preceding discussion, we must calculate $\phi(\Phi^2)$.
Since the element $r^2$ of the Weyl group of $\hat{T}_1$ can be written $w_{\alpha_*} w_{\beta_*} w_{\alpha_*} w_{\beta_*}$,
the element $\hat n = n_{\alpha_*} n_{\beta_*} n_{\alpha_*} n_{\beta_*}$ lies in the preimage of
$r^2$ in $N(\hat G, \hat T_1)$.

We first assume that $f = \hat n$, and then we will show that our result is independent of choice of lift of $r^2$.
By Lemma~\ref{tits-group}, we have $f^2 = \beta_*^\vee(-1) = \alpha_*^{\vee}(1) \beta_*^{\vee}(-1)$. Therefore, the claim follows in the case that $f = \hat{n}$.


Now suppose $f = \hat{t} \hat{n}$, for some $\hat{t} \in \hat{T}_1$.  Then $$\phi(\Phi^2) = \hat{t} \hat{n} \hat{t} \hat{n} = \hat{t} \hat{n} \hat{t} \hat{n}^{-1} \hat{n}^2.$$  We need to write $\hat{t} \hat{n} \hat{t} \hat{n}^{-1} \hat{n}^2$ in the form $\alpha_*^{\vee}(a) \beta_*^{\vee}(b)$ for some $a,b \in \mathbb{C}^*$.

We need to calculate $\hat{n} \hat{t} \hat{n}^{-1} = w_{\alpha_*} w_{\beta_*} w_{\alpha_*} w_{\beta_*}(\hat{t})$. Note that $\hat w$ acts by $(a,b) \mapsto (-a,-b)$ on $X_*(T)$.  It therefore must act by $(w,z) \mapsto (1/w,1/z)$ on $\hat T_1$.  It follows that $\hat t\hat n\hat t\hat n^{-1} = 1$.  In particular, $\phi(\Phi^2) = \alpha_*^{\vee}(1) \beta_*^{\vee}(-1)$, so we have our result.
\qed

Let $\phi$ be a TRSELP of type $r$.  We first note that this type of TRSELP is associated to the torus $T_2$.
More precisely, the Weyl group element associated to $\phi$ acts on the module $X^*(T_2) = X_*(\hat{T}_2)$ by $(a,b) \mapsto (-b,a)$, which is exactly the action that gives rise to $T_2$.  We have identified $T_2(E_2)$ with $E_2^* \times E_2^*$ in section \ref{quartictorus} via the roots $\alpha$ and $\alpha + \beta$ of $PGSp(4)$.  More precisely, our identification was

$$\varphi : T_2(E_2) \xrightarrow{\sim} E_2^* \times E_2^*$$
$$(a,b,c,d)E_2^* \mapsto (a/b, a/c) = (\alpha((a,b,c,d)E_2^*), (\alpha + \beta)((a,b,c,d)E_2^*)).$$

This gives us an identification of $\hat{T}_2$ with $\mathbb{C}^* \times \mathbb{C}^*$ that is compatible with duality of tori.  Namely, we use the images of the roots $\alpha$ and $\alpha + \beta$ of $PGSp(4)$ under the isomorphism $X^*(T_2) \xrightarrow{\sim} X_*(\hat{T}_2)$ to decompose $\hat{T}_2$.  Explicitly, the map $\mathbb{C}^* \times \mathbb{C}^* \xrightarrow{\sim} \hat{T}_2$ is given by
$$(w,z) \mapsto \alpha_*^{\vee}(w)(\alpha_*^{\vee} + \beta_*^{\vee})(z) = \alpha_*^{\vee}(wz) \beta_*^{\vee}(z).$$
Moreover, we have that $\Gal (E_2/F)$ acts on $\hat{T}_2 = \C^*\times \C^*$ via $\sigma(w,z) = (1/z, w)$.  By the theory in section \ref{groupsoftypeL}, $\phi$ naturally gives rise to a character $\chi$ of $T_2(E_2)_{\Gamma_2}$.  We wish to explicitly compute $\chi|_{\hat{H}^{-1}(\Gamma_2, T_2(E_2))}$.

Recall the exact sequence $$1 \rightarrow \hat{H}^{-1}(\Gamma_2, T_2(E_2)) \rightarrow T_2(E_2)_{\Gamma_2} \rightarrow T_2(F) \rightarrow 1$$
(since $\hat{H}^0(\Gamma_2, T_2(E_2)) = 1$).  Recall that $\phi|_{W_{E_2}}$ has image in $\hat{T}_2$.  Therefore, since the abelianization of $W_{E_2}$ is $E_2^*$, $\phi|_{W_{E_2}}$ factors through a map $\phi' : E_2^* \rightarrow \hat{T}_2$.  Recall our identification of $\hat{T}_2$ with $\mathbb{C}^* \times \mathbb{C}^*$ via $\alpha_*^{\vee}$ and $\alpha_*^{\vee} + \beta_*^{\vee}$.  We therefore get a canonical map
$$E_2^* \rightarrow \mathbb{C}^* \times \mathbb{C}^*$$
$$x \mapsto (\chi_1(x), \chi_2(x))$$
The local Langlands correspondence for tori says that  $\chi(w,z) = \chi_1(w) \chi_2(z)$ for all $w, z \in E_2^*$, where we view $\chi$ as a character of $E_2^* \times E_2^*$ via the identification of $T_2(E_2)$ with $E_2^* \times E_2^*$. Under this identification of $T_2(E_2)$ with $E_2^* \times E_2^*$, the element $(\varpi,1)$ represents the nontrivial class in $\hat{H}^{-1}(\Gamma_2, T_2(E_2))$ (see lemma \ref{T_2-H-1}).  Thus, we wish to compute
$\chi(1, \varpi)$.

Recall that the Artin map $W_{E_2} \rightarrow E_2^*$ sends $\Phi^4$ to $\varpi$. Write $\phi(\Phi^4) = \phi'(\varpi)$ as $\alpha_*^{\vee}(a) (\alpha_*^{\vee} + \beta_*^{\vee})(b)$, for some $a,b \in \mathbb{C}^*$.  The above discussion allows us to conclude that $\chi(1, \varpi) = a$.

\begin{proposition}
\label{prop:trselp2}
We have $\chi(\varpi, 1) = -1$.
Therefore, $\chi$ is a genuine character of the two-fold cover $\widetilde{T_2(F)} = T_2(E_2)_{\Gamma_2}$ of $T_2(F)$.
\end{proposition}

\proof
By the preceding discussion, we must calculate $\phi(\Phi^4)$.
The element $\hat m = n_{\alpha_*} n_{\beta_*}$ lies in the preimage of
$r$ in $N(\hat G, \hat T_2)$.  We first assume that $f = \hat m$, and then we will show that our result is independent of choice of lift of $r$.

Note that $\hat m^4 = n_{\alpha_*} n_{\beta_*} n_{\alpha_*} n_{\beta_*} n_{\alpha_*} n_{\beta_*} n_{\alpha_*} n_{\beta_*} $.
This is precisely the expression that was calculated in the proof of Lemma~\ref{tits-group} to equal $\beta^\vee(-1)$.

Therefore, in the case that $f = \hat{m}$, we have our result, since $\beta_*^{\vee}(-1) = \alpha_*^{\vee}(-1) (\alpha_*^{\vee}+\beta_*^{\vee})(-1)$, so that $\chi(\varpi,1)  = -1$.

Now suppose $f = \hat{t} \hat{m}$, for some $\hat{t} \in \hat{T}_2$.  Then $$\phi(\Phi^4) = (\hat{t} \hat{m})^4 = \hat{t} \hat{m} \hat{t} \hat{m}^{-1} \hat{m}^2 \hat{t} \hat{m}^{-2} \hat{m}^3 \hat{t} \hat{m}^{-3} \hat{m}^4.$$
Note that $\hat w$ acts by $(a,b) \mapsto (-b,a)$ on $X_*(T)$.  It therefore must act via $(w,z) \mapsto (1/z,w)$ on $\hat T_2$.  It follows that if $\hat{t}= (w,z)$, then $ \hat{t} \hat{m} \hat{t} \hat{m}^{-1} \hat{m}^2 \hat{t} \hat{m}^{-2} \hat{m}^3 \hat{t} \hat{m}^{-3} = (w,z)(1/z,w)(1/w,1/z)(z,1/w) = 1$.  Therefore, $\phi(\Phi^4) = \hat{m}^4 = \alpha_*^{\vee}(-1) (\alpha_*^{\vee} + \beta_*^{\vee})(-1)$, so we finally have our result.
\qed

To summarize, we have calculated in each case the restriction of $\chi$ to $\hat{H}^{-1}(\Gamma_i, T_i(E_i))$ by computing powers of $\phi(\Phi)$.  On the other hand, one can see that the computation of these powers of $\phi(\Phi)$ are also precisely the data needed to compute the cohomology classes of the groups of type $L$ that arise for each Langlands parameter.  In fact, Gross has predicted a canonical link between $\chi|_{\hat{H}^{-1}}$ and the cohomology class of the group of type $L$ as follows.  Recalling the notation of section \ref{groupsoftypeL}, there is a canonical sequence of isomorphisms $$H^2(\Gamma, \hat{T}) \cong H^3(\Gamma, X^*(T)) \cong \hat{H}^{-3}(\Gamma, X_*(T))^* \cong \hat{H}^{-1}(\Gamma, T(E))^*,$$ where $\hat{H}^{-1}(\Gamma, T(E))^*$ denotes the dual of $\hat{H}^{-1}(\Gamma, T(E))$.  These isomorphisms are given by the exponential sequence $$1 \rightarrow X^*(T) \rightarrow X^*(T) \otimes \mathbb{C} \rightarrow \hat{T} \rightarrow 1,$$ the universal coefficients theorem, and Tate duality, respectively.  Therefore, if $D$ is a group of type $L$, the composition of the above isomorphisms canonically gives a homomorphism $$c_D : \hat{H}^{-1}(\Gamma, T(E)) \rightarrow \mathbb{C}^*.$$  Gross has predicted that if $\phi : W_F \rightarrow D$ is a Langlands parameter, and if $\chi \in \widehat{T(E)_{\Gamma}}$ is the character that $\phi$ gives rise to as in section \ref{groupsoftypeL}, then the restriction of $\chi$ to $\hat{H}^{-1}(\Gamma,T(E))$ is equal to $c_D$.

\section{The relationship between the Gross construction and the DeBacker--Reeder construction}\label{grossdebackerreeder}

We consider here a TRSELP $\phi$ for any unramified connected reductive group G.  Let $T$ be as in \S\ref{groupsoftypeL}.  Let $E$, $\chi$, etc.~be as in~\S\ref{groupsoftypeL}.  Let $w$ be the Weyl group element associated to $\phi$, and set $\sigma = w \theta\in {\rm Aut}(T)$.  Let $\chi_{\phi}$ be the character of $T(F) = T^{\Phi_{\sigma}}$ that DeBacker and Reeder attach to $\phi$ (see section \ref{preliminaries}).

We have the exact sequence $$1 \rightarrow \hat{H}^{-1}(\Gamma, T(E)) \rightarrow T(E)_{\Gamma} \rightarrow T(F) \rightarrow \hat{H}^0(\Gamma, T(E)) \rightarrow 1$$ where $\Gamma = \Gal(E/F)$.  Recall that $\phi$ canonically gives rise to a character $\chi$ of $T(E)_{\Gamma}$, by local Langlands for tori (see section \ref{groupsoftypeL}).  Note that the above exact sequence restricts to an exact sequence
$$1 \rightarrow \hat{H}^{-1}(\Gamma, T(\mathfrak{o}_E)) \rightarrow T(\mathfrak{o}_E)_{\Gamma} \rightarrow T(\mathfrak{o}_F) \rightarrow \hat{H}^0(\Gamma, T(\mathfrak{o}_E)) \rightarrow 1$$
Moreover, one can show using a profinite version of Lang's theorem and various results about tori over finite fields, that since $T$ is unramified, $\hat{H}^{-1}(\Gamma, T(\mathfrak{o}_E)) = \hat{H}^0(\Gamma, T(\mathfrak{o}_E)) = 1$.  Therefore, the map $$T(\mathfrak{o}_E)_{\Gamma} \xrightarrow{N} T(\mathfrak{o}_F)$$ is an isomorphism, and we may view $\chi|_{T(\mathfrak{o}_E)_{\Gamma}}$ as a character of $T(\mathfrak{o}_F)$ via this isomorphism.  Now, since $\chi$ was obtained from $\phi|_{W_E}$ via the local Langlands correspondence, we get that $\chi_{\phi} \circ N = \chi$ on $T(\mathfrak{o}_E)_{\Gamma}$.  We have therefore proven the following lemma.

\begin{lemma}\label{grossanddebackerreedercompatibility}
The restriction to $T(\mathfrak{o}_E)_{\Gamma} \xrightarrow{\sim} T(\mathfrak{o}_F)$ of the genuine character arising from the Gross construction coincides with the character of $T(\mathfrak{o}_F)$ that is constructed from $\phi$ via the construction of DeBacker--Reeder construction.
\end{lemma}

Let $T(F)_{0,s}$ denote the set of strongly regular topologically semisimple elements of $T(F)$ (see \cite[\S 7]{debackerreeder}).
In the next section, we will attach a conjectural character formula to $\phi$,
and we will show that it agrees on $T(F)_{0,s}$ with the character of a unique
 depth-zero supercuspidal representation attached to $\phi$.  This character is calculated
 in~\cite{debackerreeder}, and we now recall it in a simple case.

Let $T(F)_{0,s}$ denote the set of strongly regular topologically semisimple elements of $T(F)$.  Recall that $W = N(G(F^u), T(F^u))^{\Gamma_u} / T(F)$.

\begin{definition}{\cite[\S 11]{debackerreeder}}
Suppose $\gamma \in T(F)$ is strongly regular and topologically semisimple, and let $w \in N(G(F^u), T(F)) / T(F^u)$.  Set
$${}^w \mathcal{R}_{\phi}(\gamma) = \epsilon(G, T) \displaystyle\sum_{n \in W} {}^w \chi_{\phi}(n^{-1} \gamma n),$$ where ${}^w \chi_{\phi}$ denotes the $w$-conjugate of $\chi_{\phi}$, and where $\epsilon(G,T)$ is as in \cite[\S4.4]{debackerreeder}.
\end{definition}

We set $\mathcal{R}_{\phi} = {}^1 \mathcal{R}_{\phi}$.  Associated to $\phi$, Debacker and Reeder (see \cite{debackerreeder}) have attached a collection of depth zero supercuspidal representations (i.e. an $L$-packet) indexed by $N(G(F^u), T(F)) / T(F^u)$.  We denote the depth zero supercuspidal representation associated to $w \in N(G(F^u), T(F)) / T(F^u)$ by ${}^w \pi$.  Let $\theta_{{}^w \pi}$ denote the character of ${}^w \pi$.

\begin{theorem}{\cite[\S 11]{debackerreeder}}
$\theta_{{}^w \pi}(\gamma) = {}^w \mathcal{R}_{\phi}(\gamma) \ \forall \gamma \in T(F)_{0,s}$.
\end{theorem}

\begin{remark}\label{pgsp4}
We remark that for the group that we are interested in, $PGSp(4,F)$, it is always the case that if $T$ is unramified elliptic, then $T(\mathfrak{o}_F) \cong T(F)$.  In particular, $\hat{H}^0(\Gamma, T(E)) = 1$ and the exact sequence $$1 \rightarrow \hat{H}^{-1}(\Gamma, T(E)) \rightarrow T(E)_{\Gamma} \rightarrow T(F) \rightarrow 1$$ splits.
\end{remark}

\section{Character formulas}\label{characterformulas}

\subsection{Calculating the stable character}\label{charactercalculation}
In this section, we define our character formula and prove that it agrees with the character of a unique depth-zero supercuspidal representation on $T(F)_{0,s}$ for an appropriate torus $T$.  We first recall our identification
of $T_i(E_i)$ with the group $\{ (a,b,c,d)\in E_i^*:ad=bc\}$.
Under this identifiction, we choose the following system of positive roots.
For $t = (a,b,c,d)\in T_i(E_i)$, set
\begin{equation*}
\alpha(t) = a / b,\quad \beta(t) = b / c,\quad (\alpha + \beta)(t) = a / c,\quad (2 \alpha + \beta)(t) = a / d,
\end{equation*}

In section \ref{rhoroots}, we described the groups $T_i(F)$ of $F$-rational points via the isomorphisms
\begin{eqnarray*}
\phi : T_1(E_1)  &\xrightarrow{\sim}&  E_1^* \times E_1^*\\(a,b,c,d)E_1^*  &\mapsto&  (a/b, b/c)\\ \\
\varphi : T_2(E_2) &\xrightarrow{\sim}& E_2^* \times E_2^*\\
(a,b,c,d)E_2^* &\mapsto&  (a/b, a/c) .
\end{eqnarray*}
Under these identifications, the two-fold covers $\widetilde{T_i(F)}$ arising from Langlands parameters via the theory of groups of type L are
$$\widetilde{T_1(F)} = (E_1^*\times E_1^*)/(E_1^*\times N_{E_1/F}(E_1^*)),\qquad
\widetilde{T_2(F)} = E_2^*/ N_{E_2/E_1}(E_2^*).$$
Moreover, we have identifications
$$T_1(F) = \ker(N_{E_1/F})\times \ker(N_{E_1/F}),\qquad
T_2(F) = \ker(N_{E_2/E_1}).$$
The covering norm maps are
\begin{eqnarray*}
N : \widetilde{T_1(F)} &\rightarrow & T_1(F)\\
([w_1], [w_2]) &\mapsto & (w_1 / \overline{w_1}, w_2 / \overline{w_2})\\ \\
N : \widetilde{T_2(F)} &\rightarrow & T_2(F)\\
\mbox{$[w]$} & \mapsto & w / \overline{w}
\end{eqnarray*}

We will compute our character formula using these realizations of $\widetilde{T_i(F)}$.
In particular, we will need to compute root values on elements in these renditions of $T_i(F)$.  This is done, by pulling back via $\phi^{-1}$ and $\varphi^{-1}$.  In particular, if $z \in T_1(F) = \ker(N_{E_1/F}) \times \ker(N_{E_1/F})$ and $\delta$ is a positive root as above, then $\delta(z)$ is defined by $\delta(\phi^{-1}(z))$.  If $z \in T_2(F) = \ker(N_{E_2/E_1})$, then $\delta(z) = \delta(\varphi^{-1}(z))$, where we are identifying $T_2(F)$ with $\ker(N_{E_2/E_1})$ as in the proof of Lemma \ref{other}.

In particular, it is readily computed that the values of the positive roots on $T_1(F)$ are
\begin{equation*}
\alpha((z_1,z_2)) = z_1,\quad \beta((z_1,z_2)) = z_2,\quad (\alpha + \beta)((z_1,z_2)) = z_1z_2,\quad (2 \alpha + \beta)((z_1,z_2)) = z_1^2z_2,
\end{equation*}
while those of the positive roots on $T_2(F)$ are
\begin{equation*}
\alpha(z) = z,\quad \beta(z) = \tau(z)/z,\quad (\alpha + \beta)(z) = \tau(z),\quad (2 \alpha + \beta)(z) = z\tau(z).
\end{equation*}

We now define our Weyl denominator.  Recall from section \ref{realgroups} that in the real case, the Weyl denominator was given by $\Delta^0(t, \Delta^+) \rho(\tilde{t})$, where $t \in  T(\mathbb{R})$ and $\tilde{t}$
was a lift of $t$ to $T(\mathbb{R})_{\rho}$.  Our Weyl denominator will be a $p$-adic analogue of $\Delta^0(t, \Delta^+) \rho(\tilde{t})$, as we now explain.

For a regular semisimple element $t\in T_i(F)$, note that $\Delta^0(t, \Delta^+)$ takes values in $E_i^*$, not $\mathbb C^*$.  Hence we must compose $\Delta^0$ with an appropriate $\mathbb{C}^*$-valued character $\eta$ which we define later.
For the $p$-adic version of $\rho(\tilde{t})$, we wish to define a complex-valued function on $\widetilde{T_i(F)}$ that will act as a ``square root of $\eta\circ (2 \rho)$'', in analogy with the case of real groups.
More precisely, we will define a function $\eta_{\rho}$ on $T_i(E)_{\Gamma}$ that is a canonical square root of $\eta\circ (2 \rho)\circ N$, which we will later compute in the separate case $T_1$ and $T_2$.

Before we define our character formula, we need to verify that the action of a certain Weyl group on
$T_i(E_i)$ decends to one on $T_i(E_i)_{\Gamma_i}$ and $\widetilde{T_i(F)}$.
\begin{lemma}\label{weylgroups}
Let $G$ be a connected reductive $F$-group and let $T$ be a maximal $F$-torus of $G$.  Let $E$ be the
splitting field of $T$ as in~\S\ref{groupsoftypeL}, and set $\Gamma = Gal(E/F)$.
\begin{enumerate}
\item $N(G(E), T(F)) / T(E) \cong (N(G,T)/T)(F)$.
\item The standard action of $N(G(E),T(E)) / T(E)$ on $T(E)$ determines well-defined actions of
$N(G(E), T(E))^{\Gamma} / T(F)$ and $(N(G,T)/T)(F)$ on $T(E)$,
which factor naturally to actions on $T(E)_{\Gamma}$.
\end{enumerate}
\end{lemma}

\proof
We first prove (1).  Let $nT(E)\in N(G(E), T(F)) / T(E)$.  Let $t\in T(F)$.  Then $ntn^{-1}\in T(F)$, so for $\gamma\in\Gamma$,
$\gamma(n)t\gamma(n)^{-1} =  \gamma (ntn^{-1}) = ntn^{-1}$.  Thus $n^{-1}\gamma(n)$ centralizes
$T(F)$, hence centralizes $T(E)$.  Thus $n^{-1}\gamma(n)\in T(E)$ so $\gamma(n T(E)) = \gamma (n) T(E) = n T(E)$.
Thus the elements of $N(G(E), T(F)) / T(E)$ are $\Gamma$-fixed, so $N(G(E), T(F)) / T(E) \subset (N(G,T)/T)(F)$.

Conversely, since $T$ splits over $E$, any coset in $(N(G,T)/T)(F)$ has a representative in the normalizer
$N(G(E), T(E))$.
Then $\gamma(n) = nt'$ for some $t'\in T(E)$ since $nT$ is $F$-rational.  Let $t\in T(F)$.
Then for $\gamma\in\Gamma$, $\gamma(ntn^{-1}) = \gamma(n) t \gamma(n)^{-1} = (nt')t(nt')^{-1} = ntn^{-1}$.
Thus $ntn^{-1}\in T(F)$, so $n\in N(G(E),T(F))$, as desired.

We now prove (2).
The group $N(G(E), T(E))^{\Gamma} / T(F)$ embeds naturally in $(N(G,T)/T)(F)$.  Thus it suffices to
to prove the statement for $(N(G,T)/T)(F) = N(G(E), T(F)) / T(E)$.  Let $n \in N(G(E),T(F))$.  Then since $n$
normalizes $T(F)$, $n$ also normalizes the centralizer of $T(F)$ in $G(E)$, which is  $T(E)$.  Therefore,
$N(G(E), T(F)) / T(E)$ acts on $T(E)$.

We now have to show that $N(G(E), T(F)) / T(E)$
sends the augmentation ideal $\{ t\sigma(t)^{-1} : t \in T(E), \sigma \in \Gamma \}$ to itself.  Let $t \in T(E), \sigma \in \Gamma$ and let $w \in  (N(G,T)/T)(F)$.  Then $$w \frac{t}{\sigma(t)} = \frac{wt}{w \sigma(t)} = \frac{wt}{\sigma (\sigma^{-1}(w(\sigma(t))))} = \frac{wt}{\sigma ({}^{\sigma^{-1}} w) (t)} = \frac{wt}{\sigma(w(t))},$$ the last equality coming from the fact that $w \in (N(G,T)/T)(F)$, so $w$ is fixed by Galois.
\qed

\begin{remark}
When $G = {\rm PGSp}(4)$, it is easily seen that $(N(G,T_i)/T_i)(F)$ stabilizes $\ker \chi\subset T_i(E_i)$, hence acts on $\widetilde{T_i(F)}$.
\end{remark}

Let $\phi_1$ be a TRSELP of type $r^2$ and $\phi_2$ a TRSELP of type $r$ (see Definition \ref{trselptype}). We now define a character formula associated to a Langlands parameter $\phi_i$.  Let $\chi_i \in \widehat{T_i(E_i)_{\Gamma}}$ be the character that is attached to $\phi_i$ via the theory of groups of type $L$.  Write
$$W_i = N(G(E_i), T_i(E_i))^{\Gamma_i} / T_i(F)$$ and note that $W_i$ embeds naturally in $(N(G,T_i)/T_i)(F)$.

\begin{definition}
Let $\gamma$ be a regular semisimple element of $T_i(F)$.  For $w\in (N(G,T_i)/T_i)(F)$, define
${}^w \Theta_{\chi_i}(\gamma)$
to be
$$\frac{ \displaystyle\sum_{n \in W_i} n_* w_* \chi_i(\tilde{\gamma}) }{\eta(\Delta^0(\gamma, \Delta^+)) \eta_{\rho}(\tilde{\gamma})}$$
where $\tilde{\gamma}$ is any element of $T_i(E_i)_{\Gamma_i}$ such that $N(\tilde{\gamma}) = \gamma$.
Here $s_* \chi_i(\tilde{\gamma}) := \chi_i(s^{-1} \tilde{\gamma} s)$ for $s\in (N(G,T_i)/T_i)(F)$.
We also set $\Theta_{\chi_i}(\gamma) = {}^1 \Theta_{\chi_i}(\gamma)$
\end{definition}

A priori, this expression depends on the particular choice of $\tilde\gamma$, not just $\gamma$.  However, we will verify in sections \ref{T1character} and \ref{T2character} that it doesn't depend on the choice.

We will sometimes denote $\eta(\Delta^0(\gamma, \Delta^+)) \eta_{\rho}(\tilde{\gamma})$ by
$D(\tilde{\gamma})$.

The next section concerns the character $\eta_{\rho}$ appearing in our Weyl denominator.  In the two subsequent sections, we construct $\eta_{\rho}$ for $T_1$ and $T_2$, and then we prove that $\Theta_{\chi_i}$ agrees with $\mathcal{R}_{\phi_i}$ on $T_i(F)_{0,s}$.

\subsubsection{The p-adic Weyl denominator}\label{weyldenominator}
If $i=1$, let $\eta$ be the quadratic character $\omega_{E_2/E_1}$, while if
$i=2$, let $\eta$ be either of the two characters of $E_2$ of order $4$.
In~\S\ref{T1character} and~\S\ref{T2character}, for each of the two tori $T_1$ and $T_2$, we construct a genuine character $\eta_\rho$ of
$\widetilde{T_i(F)}$ such that we have a commutative diagram
$$
\begin{CD}
\widetilde{T_i(F)} @> \eta_\rho >> \mathbb{C}^*\\
@VV N V @VV s V\\
T_i(F) @>\eta \circ (2 \rho)>> \mathbb{C}^*
\end{CD}
$$
where $s : \mathbb{C}^* \rightarrow \mathbb{C}^*$ is the map $s(z) = z^2$.
This is equivalent to showing that there exists a genuine character $\eta_\rho$ of $\widetilde{T_i(F)}$
whose square is the character
$\eta\circ(2\rho)\circ N$.
The character $\eta_\rho$ is unique up to multiplication by a non-genuine character of $\widetilde{T_i(F)}$ of order dividing $2$.
However, it is easily seen that imposing the additional condition that $\eta_\rho$ be trivial on the image of
$T_i(\mathfrak o_{E_i})_\Gamma$ in $\widetilde{T_i(F)}$ determines $\eta_\rho$ uniquely.  The character $\eta_\rho$
constructed in~\S\ref{T1character} and ~\S\ref{T2character} does indeed satisfy this extra condition and is therefore the unique such character.

\subsubsection{The torus $T_1(F)$}\label{T1character}
We first consider the torus $T_1$.  Let $\eta = \omega_{E_2/E_1}$.
Let $(z_1,z_2) \in T_1(F)$.  Then
$$\eta(\Delta^0((z_1, z_2), \Delta^+))
= \eta\left(1-\frac{1}{z_1}\right)\eta\left(1-\frac{1}{z_2}\right)\eta\left(1-\frac{1}{z_1 z_2}\right)\eta\left(1-\frac{1}{z_1^2 z_2}\right).$$
Note that this equals
$$\eta\left(1-\frac{\overline{w_1}}{w_1}\right) \eta\left(1 - \frac{\overline{w_2}}{w_2}\right)
\eta\left(1 - \frac{\overline{w_1 w_2}}{w_1 w_2}\right)
\eta\left(1 - \frac{\overline{w_1^2 w_2}}{w_1^2 w_2}\right)$$
for any $([w_1],[w_2]) \in T_1(E_1)_{\Gamma}$
such that $N([w_1],[w_2]) = (z_1, z_2)$.

We now determine the character $\eta_\rho$.
To calculate $\eta_{\rho}$, we must find a canonical square root of
$\eta\circ (2\rho)\circ N$.  For $([w_1], [w_2])\in\widetilde{T_1(F)}$, we have
$$\eta((2 \rho)(N([w_1],[w_2]))) = \eta((2 \rho)((w_1 / \overline{w_1}, w_2 / \overline{w_2})))
= \eta\left( \frac{w_1}{\overline{w_1}}
\frac{w_2}{\overline{w_2}} \frac{w_1 w_2}{\overline{w_1 w_2}} \frac{w_1^2 w_2}{\overline{w_1^2 w_2}}\right) =
\eta\left(\frac{w_1^4 w_2^3}{\overline{w_1^4 w_2^3}}\right).$$
Multiplying by $\eta(\overline{w_2})/\eta(\overline{w_2})$, and noting that
$\eta$ is trivial on $w_2\overline{w_2}$, we get
$$\eta(2 \rho(N([w_1],[w_2])))
= \eta\left(\frac{w_1^4 w_2^2}{\overline{w_1^4 w_2^4}}\right) .$$
Observe that the character
$$([w_1],[w_2])\mapsto \eta\left(\frac{w_1^2 w_2}{\overline{w_1^2 w_2^2}}\right)$$
is a square root of $\eta\circ (2\rho)\circ N$.  Moreover, this character is trivial on the
image of $T_1(\mathfrak{o}_{E_1})_{\Gamma_1}$ in $\widetilde{T_1(F)}$ and nontrivial on
$([1],[\varpi]) \in \hat H^{-1}(\Gamma_1,T_1(E_1))$, hence must equal $\eta_\rho$ (see~\ref{weyldenominator}).
Therefore, we have that
\begin{eqnarray*}
D([w_1],[w_2]) &=&
\eta\left(1-\frac{\overline{w_1}}{w_1}\right) \eta\left(1 - \frac{\overline{w_2}}{w_2}\right)
\eta\left(1 - \frac{\overline{w_1 w_2}}{w_1 w_2}\right)
\eta\left(1 - \frac{\overline{w_1^2 w_2}}{w_1^2 w_2}\right)
\eta\left(\frac{w_1^2 w_2}{\overline{w_1^2 w_2^2}}\right)\\
&=& \eta\left(\frac{1}{\overline{w_1}} - \frac{1}{w_1}\right)
\eta\left(\frac{1}{\overline{w_2}} - \frac{1}{w_2}\right)
\eta\left(\frac{1}{\overline{w_1 w_2}} - \frac{1}{w_1 w_2}\right)
\eta\left(w_1^2 w_2 - \overline{w_1^2 w_2}\right) .
\end{eqnarray*}
Since $\eta$ is trivial on the norms, we may multiply this whole expression by
$\eta(w_1 w_2 \overline{w_1 w_2})^2 = 1$ to get
$$\eta(w_1 - \overline{w_1}) \eta(w_2 - \overline{w_2})
\eta(w_1 w_2 - \overline{w_1 w_2}) \eta(w_1^2 w_2 - \overline{w_1^2 w_2}).$$

\begin{proposition}\label{quadratictoruscharactercalculation}
$\Theta_{\chi_1}(\gamma) = \mathcal{R}_{\phi_1}(\gamma)$ for all $\gamma \in  T_1(F)_{0,s}$.
\end{proposition}

\proof
Let $\gamma = (z_1, z_2) \in T_1(F)$ (so that $z_i\in\ker(N_{E_1/F})$ according to the identification in
Lemma~\ref{T_1-identification}).  Assume that $\gamma$ is strongly regular topological semisimple.
The assumption
of topological semisimplicity on $\gamma$ means that the $z_i$ are roots of unity.
Moreover, the assumption of strong regularity implies that
$z_1 \neq 1, z_2 \neq 1, z_1 z_2 \neq 1, z_1^2 z_2 \neq 1$.
We will compute the terms in $\Theta_{\chi_1}$ on $T_1(E_1)_{\Gamma_1}$.  Note that $\eta_{\rho}$ pulls back canonically to a character, which we denote $\eta_{\rho}$ again, on $T_1(E_1)_{\Gamma_1}$, via the canonical projection $T_1(E_1)_{\Gamma_1} \rightarrow \widetilde{T_1(F)}$.  Let $\tilde\gamma = ([w_1], [w_2])\in T_1(E_1)_{\Gamma_1} = (E_1^*/N_{E_1/F}(E_1^*))\times (E_1^*/N_{E_1/F}(E_1^*))$
satisfy $N(\tilde\gamma) = \gamma$.

By Remark~\ref{pgsp4}, we have that $N$ gives an isomorphism
$T_1(\mathfrak o_{E_1})_{\Gamma_1} \cong T_1(F)$.  Thus, given $\gamma$ as above,
one can choose $\tilde\gamma$ to lie in $T_1(\mathfrak o_{E_1})_{\Gamma_1}$.
We suppose first that this is the case, i.e., that $w_i\in \mathfrak{o}_{E_1}$ for $i=1,2$.
Then $w_1$ and $w_2$ may be chosen to be roots of unity.
Plugging in $w_1, w_2$, we get that $w_1 / \overline{w_1} \neq 1, w_2 / \overline{w_2} \neq 1, (w_1 w_2) / (\overline{w_1 w_2}) \neq 1, (w_1^2 w_2) / (\overline{w_1^2 w_2}) \neq 1$.  Since $w_1, w_2$ are roots of unity and $\eta$ is unramified, we therefore get that $$\eta(w_1 - \overline{w_1}) \eta(w_2 - \overline{w_2}) \eta(w_1 w_2 - \overline{w_1 w_2}) \eta(w_1^2 w_2 - \overline{w_1^2 w_2}) = 1.$$
Now recall from Section \ref{grossdebackerreeder} that
$\chi_1 = \chi_{1_{\phi_1}} \circ N$ on $T_1(\mathfrak{o}_{E_1})_{\Gamma_1}$.
Therefore,
$$\Theta_{\chi_1}(\gamma) = \frac{ \displaystyle\sum_{n \in W_1} n_* \chi_1([w_1],[w_2]) }{D([w_1],[w_2])} = \displaystyle\sum_{n \in W_1} n_* \chi_1([w_1],[w_2]) = \mathcal{R}_{\phi_1}(\gamma) $$ for $([w_1],[w_2]) \in T_1(\mathfrak{o}_{E_1})_{\Gamma_1}$ mapping to $\gamma = (z_1,z_2)$ via the norm map.

Now suppose that $w_1,w_2\in E_1^*$ are arbitrary, i.e., $\tilde\gamma$ is an arbitrary element
of $T_1(E_1)_{\Gamma_1}$ such that $N(\tilde\gamma) = \gamma$.  By Remark~\ref{pgsp4},
$\tilde\gamma = \tilde\gamma_0\delta$ for some $\gamma_0\in T_1(\mathfrak{o}_{E_1})_{\Gamma_1}$ and
$\delta\in\hat{H}^{-1}(\Gamma_1, T_1(E_1))$.  It suffices to show that
$\Theta_{\chi_1} (\gamma_0\delta) = \Theta_{\chi_1} (\gamma_0)$ for all $\delta\in\hat{H}^{-1}(\Gamma_1, T_1(E_1))$ and
strongly regular, topologically semisimple $\gamma_0\in T_1(\mathfrak{o}_{E_1})_\Gamma$.

First consider $\delta = ([\varpi], [\varpi])$.
Note that
\begin{eqnarray*}
\lefteqn
{D([\varpi w_1],[\varpi w_2])}\\
&=& \eta(\varpi w_1 - \overline{\varpi w_1}) \eta(\varpi w_2 - \overline{\varpi w_2})
\eta(\varpi w_1 \varpi w_2 - \overline{\varpi w_1 \varpi w_2}) \eta\left((\varpi w_1)^2 \varpi w_2 - \overline{(\varpi w_1)^2 \varpi w_2}\right)\\
&=& \eta(\varpi)^7 D([w_1], [w_2]) = -D([w_1],[w_2]).
\end{eqnarray*}
Moreover,
$$\sum_{n \in W_1} n_* \chi_1([\varpi w_1],[\varpi w_2]) = \sum_{n \in W_1} n_* \chi_1([\varpi],[\varpi]) n_* \chi_1([w_1],[w_2]).$$
To simplify this sum, we need to compute the Weyl group action on $T_1(E_1) = E_1^* \times E_1^*$.
Let $w_{\alpha}, w_{\beta}$ be the reflections in the Weyl group $W(G,T_1)$ of $PGSp(4)$ with respect to $T_1$
corresponding to the simple roots $\alpha$ and $\beta$.  Using the isomorphism $\phi$ as in Lemma~\ref{T_1-identification}, we have
$w_{\gamma}(w,z) := \phi(w_{\gamma}(\phi^{-1}(w,z)))$ for all $(w,z) \in E_1^* \times E_1^*$ and any root $\gamma$.
The action of $w_{\alpha}$ and $w_{\beta}$ on $T_1(E_1)$ is given by
$$w_{\alpha}((x_1, x_2, x_3, x_4)E_1^*) = (x_2,x_1,x_4,x_3)E_1^*\qquad w_{\beta}((x_1,x_2,x_3,x_4)E_1^*) = (x_1,x_3,x_2,x_4)E_1^*.$$
Thus we have
\begin{eqnarray*}
w_{\alpha}(w,z) &=& \phi(w_{\alpha}(w,1,1/z,1/(wz))) = \phi((1,w,1/(wz),1/z)) = (1/w, w^2 z)\\
w_{\beta}(w,z) &=& \phi(w_{\beta}(w,1,1/z,1/(wz))) = \phi((w,1/z,1,1/(wz))) = (wz, 1/z).
\end{eqnarray*}
It is now easily seen that $w_\alpha,w_\beta$ commute with the action of $\Gal(E_1/F)$, hence lie in $(N(G,T_1)/T_1)(F)$.
Thus $(N(G,T_1)/T_1)(F) = W(G,T_1)$.

Using the above formulas, we get that $w_{\alpha}([\varpi], [\varpi]) = ([1/\varpi], [\varpi^3])$.
But $([1/\varpi], [\varpi^3]) = ([\varpi], [\varpi]) \in E_1 / N_{E_1/F}(E_1^*) \times E_1 / N_{E_1/F}(E_1^*)$
since $\varpi^2 \in N_{E_1/F}(E_1^*)$.  Therefore, $w_{\alpha}([\varpi], [\varpi]) = ([\varpi], [\varpi])$.
Moreover, $w_{\alpha}([1], [\varpi]) = ([1], [\varpi])$.
Similarly, we have $w_{\beta}([\varpi], [\varpi]) = ([1], [\varpi])$ and, since $w_\beta$ has
order $2$, $w_{\beta}([1], [\varpi]) = ([\varpi], [\varpi])$.
It follows that we must have $w_{\alpha}([\varpi],[1]) = w_{\beta}([\varpi],[1]) = ([\varpi],[1])$.

Since the $(N(G,T_1)/T_1)(F)$ is generated by $w_{\alpha}, w_{\beta}$, since
$\chi_1([\varpi], [\varpi]) = \chi_1([1], [\varpi]) = -1$, and since $W_1 \subset (N(G,T_1)/T_1)(F)$, we have
$$\sum_{n \in W_1} n_* \chi_1([\varpi],[\varpi]) n_* \chi_1([w_1],[w_2]) = - \sum_{n \in W_1} n_* \chi_1([w_1],[w_2]).$$
Therefore, $$\frac{ \displaystyle\sum_{n \in W_1} n_* \chi_1([\varpi w_1],[\varpi w_2]) }{D([\varpi w_1],[\varpi w_2])} = \frac{ \displaystyle\sum_{n \in W_1} n_* \chi_1([w_1],[w_2]) }{D([w_1],[w_2])}.$$

We now consider the element $([1], [\varpi])$.  It is easy to see that $D([w_1], [\varpi w_2]) = \eta(\varpi)^3 D([w_1],[w_2])$.
Since $\eta(\varpi) = -1$, we have $D([w_1],[\varpi w_2]) = - D([w_1],[w_2])$.  Also, the above discussion of the action
of $(N(G,T_1)/T_1)(F)$ shows that
$$\displaystyle\sum_{n \in W_1} n_* \chi_1([1],[\varpi]) n_* \chi_1([w_1],[w_2]) =
- \displaystyle\sum_{n \in W_1}  n_* \chi_1([w_1],[w_2]),$$ which handles the case $\delta = ([1], [\varpi])$.

Now consider $([\varpi], [1])$.  It is easy to see that $D([\varpi w_1],[ w_2]) = \eta(\varpi)^4 D([w_1], [w_2]) = D([w_1],[w_2])$.  Moreover,
$$\sum_{n \in W_1} n_* \chi_1([\varpi w_1],[w_2]) = \sum_{n \in W_1} n_* \chi_1([\varpi],[1]) n_* \chi_1([w_1],[w_2]) .$$
Since $w_\alpha, w_\beta$ both fix $([\varpi],[1])$, and since $\chi_1([\varpi],[1]) = 1$, we have
$$\frac{ \displaystyle\sum_{n \in W_1} n_* \chi_1([\varpi w_1],[w_2]) }{D([\varpi w_1],[w_2])} = \frac{ \displaystyle\sum_{n \in W_1} n_* \chi_1([w_1],[w_2]) }{D([w_1],[w_2])}.$$
\qed

\subsubsection{The torus $T_2(F)$}\label{T2character}
Let $\eta$ be an unramified character of $E_2^*$ of order $4$.
Let $z$ be an element of $\ker (N_{E_2/E_1})$, which we identify with $T_2(F)$ as in
Lemma~\ref{other}.  Then
\begin{equation}
\label{eq:Delta0}
\eta(\Delta^0(z, \Delta^+))
= \eta\left(1 - \frac{1}{z}\right)\eta\left(1-\frac{z}{\tau(z)}\right)\eta
\left(1-\frac{1}{\tau(z)}\right)\eta\left(1-\frac{1}{z \tau(z)}\right).
\end{equation}
Let $[w]\in E_2^*/N_{E_2/E_1}(E_2^*) = \widetilde{T_2(F)} = T_2(E_2)_{\Gamma_2}$ satisfy $N(w)=w/\tau^2 (w) = z$.
Then~\eqref{eq:Delta0} equals
$$ \eta\left(1 - \frac{\tau^2(w)}{w}\right)\eta\left(1-\frac{w\tau^3(w)}{\tau(w)\tau^2(w)}\right)\eta
\left(1-\frac{\tau^3(w)}{\tau(w)}\right)\eta\left(1-\frac{\tau^2(w)\tau^3(w)}{w \tau(w)}\right).$$

We now wish to define the analogous function $\eta_{\rho}$ on $T_2(E_2)_{\Gamma_2}$ as we did for
$T_1(E_1)_{\Gamma_1}$.  Again, we first compute $\eta\circ (2 \rho)\circ N$.  For
$[w]\in  T_2(E)_\Gamma$,
let $z = N([w]) = w/\tau^2(w)$.  Then $\eta((2 \rho)(z)) = \eta( z\tau(z)^3)$
so
$$\eta((2 \rho)(N([w])) = \eta \left(\frac{w\tau(w)^3}{\tau^2(w)\tau^3(w)^3}\right).$$
If we multiply this by $\eta(w \tau(w))/\eta(w \tau(w))$, we get
$$\eta((2 \rho)(N([w]))
= \eta\left(\frac{w^2 \tau(w)^4}{w \tau(w) \tau^2(w) \tau^3(w) \tau^3(w)^2}\right) =
\eta\left(\frac{w^2 \tau(w)^4}{\tau^3(w)^2}\right)$$
since $\eta$ is trivial on $N_{E_2/F}(E_2^*)$.
Thus the character
$$w\mapsto \eta\left(\frac{w \tau(w)^2}{\tau^3(w)}\right)$$
is a square root of $\eta\circ (2 \rho)\circ N$.
This character is trivial on the
image of $T_2(\mathfrak{o}_{E_2})_{\Gamma_2}$ and nontrivial on
$\hat H^{-1}(\Gamma_2,T_2(E_2))$, hence must equal $\eta_\rho$ (see~\ref{weyldenominator}).

Therefore, we get
\begin{eqnarray}
\label{eq:Delta}
D([w]) &=& \eta(\Delta^0(N([w]), \Delta^+)) \eta_{\rho}([w])\nonumber \\
&=&
 \eta\left(1 - \frac{\tau^2(w)}{w}\right)\eta\left(1-\frac{w\tau^3(w)}{\tau(w)\tau^2(w)}\right)\eta
\left(1-\frac{\tau^3(w)}{\tau(w)}\right)\eta\left(1-\frac{\tau^2(w)\tau^3(w)}{w \tau(w)}\right)
\eta\left(\frac{w \tau(w)^2}{\tau^3(w)}\right)\nonumber\\
&=& \eta\left(w - \tau^2(w)\right) \eta\left(\tau(w) - \frac{w \tau^3(w)}{\tau^2(w)}\right)
\eta\left(\tau(w) - \tau^3(w)\right) \eta\left(\frac{1}{\tau^3(w)} - \frac{\tau^2(w)}{w \tau(w)}\right)
\end{eqnarray}
Since $\eta$ is trivial on $N_{E_2/F}(E_2^*)$, we have $\eta(w \tau(w) \tau^2(w) \tau^3(w))=1$.
Thus, multiplying \eqref{eq:Delta} by $\eta(w \tau(w) \tau^2(w) \tau^3(w))$, we find that
\begin{equation}
\label{eq:T2denominator}
D([w]) = \eta(w - \tau^2(w)) \eta(\tau(w)\tau^2(w) - w \tau^3(w)) \eta(\tau(w) - \tau^3(w)) \eta(w \tau(w) - \tau^2(w) \tau^3(w)).
\end{equation}

\begin{proposition}\label{quartictoruscharactercalculation}
$\Theta_{\chi_2}(\gamma) = \mathcal{R}_{\phi_2}(\gamma) \ \forall \gamma \in  T_2(F)_{0,s}$.
\end{proposition}

\proof
Let $\gamma \in T_2(F)$ correspond to the element $z\in \ker( N_{E_2/E_1})$ via our identification of $T_2(E_2)$ in section \ref{quartictorus}.  Suppose that $\gamma$ is strongly regular and topological semisimple.  Then $z$ is a root of unity,
and moreover, $z \neq 1$, $z \neq \tau(z)$, and $z\neq \tau(z)^{-1}$.  Let $\tilde\gamma\in T_2(E_2)_{\Gamma_2}$ satisfy
$N(\tilde\gamma) = \gamma$.  Then $\tilde\gamma$ corresponds to an element $[w]\in E_2^*/N_{E_2/E_1}(E_2^*)$ (as
in Lemma~\ref{other}) such that $w/ \tau^2(w) = z$.  The strong regularity of $\gamma$ implies that $w\neq\tau^2(w)$,
$w\tau^3(w)\neq \tau(w)\tau^2(w)$, and $w\tau(w)\neq\tau^2(w)\tau^3(w)$.

As in the case of the torus $T_1$, by Remark~\ref{pgsp4}, $N$ gives an isomorphism
$T_2(\mathfrak o_{E_2})_{\Gamma_2} \cong T_2(F)$.  Thus, given $\gamma$ as above,
one can choose $\tilde\gamma$ to lie in $T_2(\mathfrak o_{E_2})_{\Gamma_2}$.
Suppose first that this is the case, i.e., that $w\in \mathfrak{o}_{E_2}$.
Then $w$ may be chosen to be a root of unity.  It follows from this and \eqref{eq:T2denominator} that $D([w]) = 1$.
Since $\chi_2 = \chi_{2_{\phi_2}} \circ N$ on $T_2(\mathfrak{o}_{E_2})_{\Gamma_2}$ (see~\S\ref{grossdebackerreeder}), we obtain
$$\Theta_{\chi_2}(\gamma) = \frac{ \displaystyle\sum_{n \in W_2} n_* \chi_2([w]) }{D([w])} = \displaystyle\sum_{n \in W_2} n_* \chi_2([w]) = \mathcal{R}_{\phi_2}(\gamma). $$

Now suppose that $w\in E_2^*$ is arbitrary, i.e., $\tilde\gamma$ is an arbitrary element
of $T(E_2)_\Gamma$ such that $N(\tilde\gamma) = \gamma$.  As in~\S\ref{T1character}, it suffices to
show that $\Theta_{\chi_2} (\gamma_0\delta) = \Theta_{\chi_2} (\gamma_0)$ for all
$\delta\in\hat{H}^{-1}(\Gamma_2, T_2(E_2))$ and
strongly regular, topologically semisimple $\gamma_0\in T_2(\mathfrak{o}_{E_2})_{\Gamma_2}$.
From Lemma~\ref{T_2-H-1}, it follows that the nontrivial element $\delta$ of $\hat{H}^{-1}(\Gamma_2, T_2(E_2))
\subset T_2(E_2)_{\Gamma_2}$ is $[\varpi]$.

First note that
\begin{eqnarray*}
\lefteqn{D([\varpi w])}\\
&=& \eta(\varpi w - \tau^2(\varpi w)) \eta(\tau(\varpi w)\tau^2(\varpi w) - \varpi w \tau^3(\varpi w))
\eta(\tau(\varpi w) - \tau^3(\varpi w))\\
& & \cdot\ \eta(\varpi w \tau(\varpi w) - \tau^2(\varpi w) \tau^3(\varpi w))\\
&=& \eta(\varpi)^6 \eta(w - \tau^2(w)) \eta(\tau(w)\tau^2(w) - w \tau^3(w))
\eta(\tau(w) - \tau^3(w)) \eta(w \tau(w) - \tau^2(w) \tau^3(w))\\
&=& \eta(\varpi)^6 D([w]).
\end{eqnarray*}
This is equal to $-D([w])$ since $\eta(\varpi)$ has order $4$.

To simplify the numerator of $\Theta_{\chi_2}([\varpi w])$, we compute the Weyl group action on $T_2(E_2) = E_2^* \times E_2^*$.
Let $w_{\alpha}, w_{\beta}$ be the reflections in the Weyl group $W(G,T_2)$ of $PGSp(4)$ with respect to $T_2$
corresponding to the simple roots $\alpha$ and $\beta$.
Under the identification $\varphi$ Lemma~\ref{T2-identification}, we have
\begin{eqnarray*}
w_{\alpha}(w,z) = \varphi(w_{\alpha}(wz,z,w,1)) = \varphi(z,wz,1,w) = (1/w, z)\\
w_{\beta}(w,z) = \varphi(w_{\beta}(wz,z,w,1)) = \varphi((wz,w,z,1)) = (z,w).
\end{eqnarray*}
Recall from Lemma~\ref{T2-identification}, that $\tau (w,z) = (\tau(z)^{-1},\tau(w))$.  One easily checks that the subgroup $(N(G,T_2)/T_2)(F)$ of
elements of $W(G,T_2)$ which commute with the action of $\Gal (E_2/F)$ is precisely $\langle w_\alpha w_\beta\rangle$.
Now the element $[\varpi]$ of $E_2^*/N_{E_2/E_1}(E_2^*)$ corresponds to the coset in $T_2(E_2)_{\Gamma_2}$ represented by the element $(\varpi,1)$.  But $(w_\alpha w_\beta) (\varpi,1) = (1,\varpi)$, which corresponds to the element
$[\varpi^{-1}] = [\varpi]\in E_2^*/N_{E_2/E_1}(E_2^*)$.  If follows that $(N(G,T_2)/T_2)(F)$ acts trivially on
$\hat{H}^{-1}(\Gamma_2, T_2(E_2))$.  Therefore, since $W_2 \subset (N(G,T_2)/T_2)(F)$,
$$
\sum_{n \in W_2} n_* \chi_2([\varpi w]) = \sum_{n \in W_2} n_* \chi_2([\varpi]) n_* \chi_2([w]) =  \sum_{n \in W_2} \chi_2([\varpi]) n_* \chi_2([w]) = - \sum_{n \in W_2} n_* \chi_2([w])$$ since $\chi_2([\varpi]) = -1$.
Therefore, we get $$\frac{ \displaystyle\sum_{n \in W_2} n_* \chi_2([ \varpi w]) }{D([ \varpi w])} =
\frac{ -\displaystyle\sum_{n \in W_2} n_* \chi_2([w]) }{-D([w])} =
\frac{ \displaystyle\sum_{n \in W_2} n_* \chi_2([w]) }{D([w])}$$ and we are done.
\qed

\begin{remark}\label{choiceofpositiveroots}
We note that our definition of $\Theta_{\chi_i}$ involved a choice of positive roots.  However, after some calculation, one can see that the functions $\Theta_{\chi_i}$ do not depend on the choice of positive roots.
\end{remark}

\subsection{Uniqueness of restrictions of characters}\label{restriction}

In this section, we show that a depth-zero supercuspidal character of $PGSp(4,F)$ that comes from the torus $T_i(F)$
is uniquely determined by its restriction to $T_i(F)_{0,s}$.  We do this in the setting of a general unramified group.
From now on, let $G$ denote an unramified reductive
$F$-group.  Let $\pi$ be a depth-zero supercuspidal representation of $G(F)$ associated to an elliptic maximal
$F$-torus $T$.  We show that if the character of $\pi$ coincides on $Z(F) T(F)_{0,s}$ with that of another depth-zero
supercuspidal representation of $G(F)$, then these representations are equivalent.
We do this first when both representations arise from the same elliptic maximal $F$-torus $T$ of $G$.
Recall that $W$ denotes the Weyl group $N(G(F^u),T(F^u))^{\Gamma_u} / T(F)$.

\begin{lemma}\label{nonvanishing}
Let $\chi_1$ be a depth zero character of $T(F)$.  For $q$ sufficiently large, the function $$\sum_{n \in W} n_* \chi_1$$ is not identically zero on $Z(F) T(F)_{0,s}$.
\end{lemma}
\proof
Suppose to the contrary that $f = \sum_{n\in W}  n_* \chi_1$
vanishes on $Z(F) T(F)_{0,s}$.
The group $\mathbb T^{\Phi_\sigma}\subset\mathbb T$ is the quotient of the maximal compact subgroup
$T(\mathfrak o_F)$ of $T(F)$ by its pro-unipotent radical.  Since $\chi_1$ has depth zero, $f$ can be viewed as a function on $\mathbb T^{\Phi_\sigma}$.  Since $\chi_1$ is regular, the characters $n_*\chi_1$ are all distinct as $n$ ranges over $W$.  Thus
$\langle f, \chi_1\rangle = 1$, where $\langle\cdot ,\cdot\rangle$ denotes the standard Hermitian inner product on the space of complex-valued functions on $\mathbb T^{\Phi_\sigma}$.

Let $Y$ be the complement of the image of $T(F)_{0,s}\cap T(\mathfrak o_F)$ in $\mathbb T^{\Phi_\sigma}$.  Then as a function on $\mathbb T^{\Phi_\sigma}$, $f$ vanishes off of $Y$.  Thus
$$1= \langle f, \chi_1\rangle = \frac{1}{|\mathbb T^{\Phi_\sigma}|}\sum_{y\in Y}\sum_{n\in W} (n_* \chi_1)(y)\cdot \chi_1^{-1}(y)
\leq \frac{|Y||W|}{|\mathbb T^{\Phi_\sigma}|}.$$
It follows that
\begin{equation}
\label{eq:bound}
|Y|/|\mathbb T^{\Phi_\sigma}|\geq 1/|W|.
\end{equation}
It is easily seen that $Y$ lies in the union $\mathbb Y$ of the fixed-point subgroups $\mathbb T^w$ of $\mathbb T$ as $w$ ranges over $W\setminus \{ 1\}$.  Note that $\mathbb Y$ is a subvariety of $\mathbb T$ of dimension strictly smaller than that of $\mathbb T$ whose definition is independent of the particular residue field $\mathfrak f$.  It follows that $|Y|/|\mathbb T^{\Phi_\sigma}|\rightarrow 0$ as $q$ increases.  Hence, for $q$ sufficiently large, $|Y|/|\mathbb T^{\Phi_\sigma}| < 1/|W|$.  For such $q$, \eqref{eq:bound} implies that $f$ cannot be identically $0$ on $Z(F) T(F)_{0,s}$.
\qed

\begin{remark}\label{nonvanishing1}
In the case of $G = PGSp(4)$, one can calculate that the inequality $|Y|/|\mathbb T^{\Phi_\sigma}| <1/|W|$ holds for $T_1$ when $q > 46$ and for $T_2$ when $q > 3$.
\end{remark}

\begin{proposition}\label{samecartans}
Let $\chi_1, \chi_2$ be depth zero characters of $\widehat{T(F)}$.
Suppose $$\displaystyle\sum_{n \in W} n_* \chi_1(z \gamma_0) = \displaystyle\sum_{n \in W} n_* \chi_2(z \gamma_0)$$
for all strongly regular, topologically semisimple elements $\gamma_0 \in T(F)_{0,s}$ and for all elements $z \in Z(F)$.  Then $\chi_1 = n'_* \chi_2$ for some $n' \in W$.
\end{proposition}

\proof
Let $S := T(F) \setminus Z(F) T(F)_{0,s}$.  By Lemma~\ref{nonvanishing},
$$\displaystyle\sum_{n \in W} n_* \chi_1(\gamma)$$ is non-zero for some $\gamma \in Z(F) T(F)_{0,s}$
Since $\chi_1$ is depth zero, it follows in particular that $\chi_1|_S = n_* \chi_1|_S$ for all $n \in W$.  Therefore, the same argument as in the proof of \cite[Lemma 5.1]{spice} shows that there exists $n'' \in W$ such that $\chi_1 = n''_* \chi_2$ on $S$ (note that \cite[Lemma 5.1]{spice} applies when $W$ is of prime order.  However, the proof holds for general $W$ since in our case we have that $\chi_1|_S = n_* \chi_1|_S$ for all $n \in W$).

We now have that $$\displaystyle\sum_{n \in W} n_* \chi_1(\gamma) = \displaystyle\sum_{n \in W} n_* \chi_2(\gamma) \ \ \forall \gamma \in T(F)$$
Therefore, by linear independence of characters we get that $\chi_1 = n'_* \chi_2$ for some $n' \in W$.
\qed

We now show that the character of the above supercuspidal representation $\pi$ associated to $T$ cannot cannot
agree on $Z(F) T(F)_{0,s}$ with the character of a depth-zero supercuspidal representation associated to
a maximal torus that is not $F$-conjugate to $T$.

\begin{proposition}\label{differentcartans}
Suppose that $\varphi, \varphi'$ are TRSELPs associated to non-conjugate Weyl group elements $w, w'$.
Let $T, T'$ denote the elliptic tori associated to $w, w'$.  Then there exists $\gamma \in Z(F) T(F)_{0,s}$ such that $\mathcal{R}_\varphi(\gamma) \neq \mathcal{R}_{\varphi'}(\gamma)$.
\end{proposition}

\proof
It is consequence of remarks in \cite[\S 10.1]{debackerreeder} that $\mathcal{R}_{\varphi'}$ vanishes completely on $Z(F) T(F)_{0,s}$.  The reason is that if $\gamma_0 \in T(F)_{0,s}$, then $\gamma_0$ is not contained in any conjugate of $T'$. But we have shown in Lemma \ref{nonvanishing} that $\mathcal{R}_{\varphi}$ doesn't vanish on all of $Z(F) T(F)_{0,s}$.
\qed

There is one minor point here to resolve.  We do not claim that there are no other supercuspidal representations outside of \cite{debackerreeder} of $G(F)$ whose characters agree with $\Theta_{\chi}$ on $Z(F) T(F)_{0,s}$, even though we expect this to be true.  The reason we cannot claim this is that these characters haven't yet been completely computed in the literature.

We are ready to state our main theorem.  Let $\phi$ be a TRSELP for $PGSp(4,F)$ with associated Weyl group element $w$, and let $T(F)$ be the elliptic torus
of $PGSp(4,F)$ defined by $w$, and suppose $T$ splits over $E$, with $\Gamma := \Gal(E/F)$.  Recall that to $\phi$ we may canonically associate a character $\chi \in \widehat{T(E)_{\Gamma}}$ as in section \ref{groupsoftypeL}. Combining propositions \ref{samecartans}, \ref{differentcartans}, \ref{quadratictoruscharactercalculation} and \ref{quartictoruscharactercalculation}, we may now set $L(\chi)$ to denote the collection $\{ {}^w \sigma : w \in N(G(E), T(F)) / T(E) \}$, where ${}^w \sigma$ is the unique depth zero supercuspidal representation of $PGSp(4,F)$ considered in \cite{debackerreeder}, whose restriction to $T(F)_{0,s}$ equals ${}^w \Theta_{\chi}$.

\begin{theorem}
The composite map
\begin{equation*}
\phi  \mapsto  \chi \in \widehat{T(E)_{\Gamma}} \mapsto L(\chi)
\end{equation*}
is the tame local Langlands correspondence for $PGSp(4,F)$.
\end{theorem}


\begin{thebibliography}{9}

\bibitem{adams}
  Jeffrey Adams,
  \emph{Extensions of tori in ${\rm SL}(2)$.} Pacific J. Math. 200 (2001), no. 2, 257--271.

\bibitem{adamsvogan}
  Jeffrey Adams and David Vogan
  L-Groups, Projective Representations, and the Langlands Classification.
  Amer. Journal of Math. 113 (1991), 45-138.

\bibitem{adler}
  Jeffrey Adler,
  \emph{Refined anisotropic K-types and supercuspidal representations}, Pacific J. Math., 185 (1998), pp. 1-32.

\bibitem{adrian}
  Moshe Adrian
  \emph{On the Local Langlands Correspondences of DeBacker/Reeder and Reeder for $GL(\ell,F)$, where $\ell$ is prime}, Pacific Journal of Mathematics 255-2 (2012), 257--280.

\bibitem{adrian1}
  Moshe Adrian
  \emph{A New Construction of the Local Langlands Correspondence for $GL(n,F)$, $n$ a prime}, Ph.D. Thesis.

\bibitem{bushnellhenniart}
  C. Bushnell and G. Henniart,
  \emph{The essentially tame local Langlands correspondence, I},
    J. Amer. Math. Soc. 18 (2005), no. 3, 685--710.

\bibitem{carter}
  \emph{R. Carter},
  Finite Groups of Lie Type: Conjugacy Classes and Complex Characters, John Wiley and Sons Inc, 1993.

\bibitem{carter}
  \emph{R. Carter},
  Simple Groups of Lie Type, John Wiley and Sons Inc, 1993.

\bibitem{debackerreeder}
  Stephen DeBacker and Mark Reeder,
  \emph{Depth-zero supercuspidal $L$-packets and their stability.}
  Ann. of Math. (2) 169 (2009), no. 3, 795--901.

\bibitem{gross}
  Benedict Gross, \emph{Groups of type L}, unpublished notes.

\bibitem{grossreeder}
  Benedict Gross and Mark Reeder,
  \emph{Arithmetic invariants of discrete Langlands parameters.}  Duke Math. Journal, 154, (2010), 431-508.

\bibitem{kaletha1}
  Tasho Kaletha, \emph{Simple Wild L-packets}, preprint.

\bibitem{morris}
  Lawrence Morris,
  \emph{Some tamely ramified supercuspidal representations of symplectic groups.} Proc. London Math. Soc. (3) 63 (1991), no. 3, 519--551.

\bibitem{reeder}
  M. Reeder,
  \emph{Supercuspidal $L$-packets of positive depth and twisted Coxeter elements},
  J. Reine Angew. Math. 620 (2008), 1-33.

\bibitem{spice}
  Loren Spice,
  \emph{Supercuspidal characters of $SL_{\ell}$ over a $p$-adic field, $\ell$ a prime.} Amer. J. Math. 121 no. 1, 51–100.

\end{thebibliography}
\end{document}